\documentclass{mpaper}
\usepackage{tikz}
\usepackage{graphicx} % Required for inserting images

\input{symbols}
\usepackage[shortlabels]{enumitem}

\usepackage{todonotes} 

\title{Well begun is half done: the computation of Eulerian paths}
\author[Valenti]{Manlio Valenti}

\address[Manlio Valenti]{Department of Computer Science\\
Swansea University\\
Swan\-sea, SA1 8EN\\
UK}
\email{\href{mailto:manliovalenti@gmail.com}{manliovalenti@gmail.com}}

\date{}

\thanks{The author thanks Jeffry Hirst and Alberto Marcone for useful conversations on the topics of the paper. ChatGPT-6 was used to explore the literature, as well as discuss and critically review some of the proof strategies, and polish the presentation. All suggestions were carefully assessed, and the author takes full responsibility for the correctness of the results.}

\begin{document}
\maketitle

%    Abstract is required.
\begin{abstract}
We characterize the uniform computational strength of finding (one-way or two-way, directed or undirected) Eulerian paths in countable graphs. In the one-way case, for both directed and undirected graphs, this problem is Weihrauch-equivalent to $\mflim\times\CNatural'$, where $\CNatural'$ is the jump of closed choice on the natural numbers. This refines the proof-theoretic strength of the existence of a one-way Eulerian path given by the reverse mathematics framework. The $\CNatural'$ component comes entirely from the difficulty of choosing a starting vertex: with a valid start prescribed, the problem is equivalent to $\mflim$,  which also characterizes the degree of the two-way directed/undirected Eulerian path problem. 

We also study the complexity of orienting the edges of a graph while preserving the existence of a one-way/two-way Eulerian path. In the one-way case, this problem is not computable, and we prove that $\mflim_2\strictlyweireducible\OrientOW\strictlyweireducible\Choice{2}'\sqcap\Sort_2$ and $\OrientOW\not\weireducible\CNatural$. On the other hand, valid orientations are uniformly computable for two-way paths and for one-way paths with a prescribed valid starting vertex.
\end{abstract}

\tableofcontents

\section{Introduction}
An \textdef{Eulerian path} in a graph is a trail that traverses every edge exactly once. For infinite graphs, one distinguishes between one-way paths, indexed by $\mathbb{N}$, and two-way paths, indexed by $\mathbb{Z}$. The classical Erd\H{o}s-Gr\"unwald-Weiszfeld theorem (see \cite[Thm.\ 3.2.1]{Ore62}) provides a characterization of when a countable graph admits a one-way Eulerian path. Analogous characterizations can be obtained for one-way Eulerian paths on directed graphs \cite[Thm.\ 1]{NashWilliams66}, and two-way Eulerian paths for undirected graphs \cite[Thm.\ 3.2.2]{Ore62} and directed graphs \cite[Thm.\ 2]{NashWilliams66}.

Eulerian paths have been studied in computable graph theory and reverse mathematics. In \cite[Thm.\ 1]{Bean76}, the author observed that there are computable graphs with an Eulerian path, but no computable such path. In contrast, every highly recursive graph (i.e.\ a locally finite computable graph admitting a computable degree function $d\function{V}{\mathbb{N}}$ mapping every vertex to the number of its neighbors) with an Eulerian path has a computable such path. This result was recently generalized to moderately computable graphs (i.e.\ computable graphs with a computable degree function $d\function{V}{\mathbb{N}\cup \{\infty\}}$) in \cite{CarrascoVargas25}. In \cite{CVDRR26}, the authors study the computability of separating infinite graphs into infinite components and how the number of ends affects the complexity of finding infinite simple paths and deciding whether Eulerian paths exist.

From the point of view of reverse mathematics, it has been shown that, even for pre-Eulerian (i.e.\ satisfying the hypotheses of the Erd\H{o}s-Gr\"unwald-Weiszfeld theorem) locally finite graphs, the existence of an Eulerian path is equivalent to $\mathrm{ACA}_0$ over $\mathrm{RCA}_0$ \cite[Thm.\ 12]{GH98}. However, a closer look at the proof reveals that, while there is a computable graph $G$ such that every Eulerian path through $G$ computes $\emptyset'$, the proof of the existence of an Eulerian path through every pre-Eulerian graph uses $\emptyset''$, therefore leaving a gap. 

In this paper, we use the framework of Weihrauch reducibility to analyze the uniform computational strength of the problem of finding (one-way or two-way) Eulerian paths through countable (directed or undirected) graphs. 

\begin{definition}
	Let $\EP$ be the problem that takes as input a countable undirected graph that has a one-way Eulerian path, and outputs a sequence $\sequence{x_n}{n\in\mathbb{N}}$ of vertices forming a one-way Eulerian path, and let $\DEP$ be the analogous problem on directed graphs. We also write $\TEP$ (resp.\ $\TDEP$) for the problem that takes as input a countable undirected (resp.\ directed) graph that has a two-way Eulerian path, and outputs a two-way undirected (resp.\ directed) Eulerian path $\sequence{x_n}{n\in\mathbb{Z}}$.
\end{definition}

We assume that a graph is represented via the characteristic function of its set of vertices and edges.

In \thref{thm:EP}, we show that 
\[
\EP\weiequiv\DEP\weiequiv\mflim\times\CNatural',
\]
where $\mflim$ is the limit operator on Baire space and $\CNatural'$ is the jump of closed choice on $\mathbb{N}$ (see Section~\ref{sec:weihrauch}). This refines the proof-theoretic characterization obtained via the reverse mathematics analysis: while every computable graph $G\in \dom(\EP)$ has a $\emptyset'$-computable solution, no uniform $\emptyset'$-computable procedure produces one. We observe that the gap between $\EP$ (or $\DEP$) and $\mflim$ rests entirely on the complexity of finding the first vertex of the Eulerian path. Once such a vertex is found, the rest of the path can be obtained uniformly in the jump of the input (\thref{thm:fixed-start}). This also gives an intuitive explanation of why, in contrast to what happens for one-way paths, finding two-way (directed or undirected) Eulerian paths is simpler (\thref{thm:tep=lim} and \thref{thm:tdep=lim}):
\[ \EP_2\weiequiv\DEP_2\weiequiv\mflim. \]

The similarity of the directed and undirected constructions suggests a second question: can we computably orient an undirected graph while preserving the existence of a (one-way or two-way) Eulerian path? In Section~\ref{sec:orientations}, we study the Weihrauch degree of the problem $\OrientOW$ of computing an orientation of the edges of a graph preserving one-way Eulerian paths. We obtain  
\[
\mflim_2\strictlyweireducible\OrientOW
\strictlyweireducible\Choice{2}'\sqcap\Sort_2,
\qquad \OrientOW\not\weireducible\CNatural.
\]
%$\mflim_2 \strictlyweireducible \OrientOW \strictlyweireducible \Choice{2}' \sqcap \Sort_2$ (\thref{thm:lim_2<orient} and \thref{thm:orient<meet}). 
The $\Choice{2}'$ upper bound is obtained by showing that, for every undirected graph $G\in\dom(\EP)$, we can uniformly compute a directed graph $D$ by orienting the edges of $G$ such that either $D$ or its global reversal $D^\mathrm{op}$ is a valid solution for $\OrientOW(G)$ (\thref{thm:D_Dop}). This implies, e.g.\ that $\Choice{3}\not\weireducible \OrientOW$. The obtained results separate $\OrientOW$ from other standard benchmarks in the Weihrauch lattice. 

Once again, the complexity of $\OrientOW$ can be explained in terms of the complexity of finding the starting vertex of the path. In contrast, the problem $\OrientTW$ of orienting the edges of an undirected graph so that the resulting directed graph admits a two-way Eulerian path is uniformly computable (\thref{thm:computable-orientation-two-way}).

\section{Background and notation}

We briefly recall the main notions and fix the notation that will be used throughout the paper. 

An \textdef{undirected graph} $G$ is a pair $(V,E)$ where $V$ is the set of \textdef{vertices} and $E$ is a set of unordered pairs of vertices, called \textdef{edges}. A \textdef{directed graph} is a graph where the edges are ordered pairs. To stress the difference between unordered and ordered edges between two vertices $u$ and $v$, we use $\{u,v\}$ to denote the former and $(u,v)$ to denote the latter. In this paper, we will always assume graphs to be \textdef{simple}, i.e.\ graphs with no self-loops such that there can be at most one edge between any (ordered or unordered, depending on whether the graph is directed or not) pair of vertices. Given a graph $G$, we use $V(G)$ and $E(G)$ to denote, respectively, the set of vertices and the set of edges of $G$. We write $\deg_G(v)$ for the \textdef{degree of a vertex $v$ in $G$}, namely the number of neighbors of $v$.

A \textdef{path} or \textdef{walk} through a graph $G$ is a (finite or infinite) sequence of vertices $\sequence{v_n}{n}$ of $G$ such that, for each $n$, the edge $\{v_n,v_{n+1}\}$ (or $(v_n,v_{n+1})$ in case $G$ is directed) is in $E(G)$. Since we are working with simple graphs, a path uniquely identifies a sequence of edges. A \textdef{trail} is a path where no edge is repeated. A finite trail is called \textdef{closed} or \textdef{open} according to whether the initial and final vertices are the same or not. A \textdef{ray} is an infinite or finite and open trail with no repeated vertices. A \textdef{cycle} is a closed (non-empty) trail with at least one edge and no repeated vertices except for its endpoints. An \textdef{Eulerian path} is a trail that uses each edge and goes through each vertex. We distinguish between \textdef{one-way Eulerian paths}, namely Eulerian paths of the form $\sequence{v_n}{n\in\mathbb{N}}$ and \textdef{two-way Eulerian paths}, i.e.\ Eulerian paths of the form $\sequence{v_n}{n\in\mathbb{Z}}$. 

There are well-known results that characterize exactly the graphs that admit Eulerian paths. To state them, we first introduce the following definitions: given a graph $G$ and a set of edges $F$, we define $G-F$ as the graph obtained by removing all edges in $F$ from $G$, as well as every vertex that remains isolated after the removal of $F$. In symbols:
\begin{gather*}
	E(G-F):=E(G)\setminus F,\\
	V(G-F):= V(G) \setminus \set{v \in V(G)}{v \text{ is isolated in } (V(G), E(G-F))}.
\end{gather*}

\begin{theorem}[{Erd\H{o}s-Gr\"unwald-Weiszfeld \cite[Thm.\ 3.2.1]{Ore62}}]
	\thlabel{thm:EGW}
A countable graph $G$ has an infinite one-way undirected Eulerian path iff 
\begin{enumerate}
    \item $G$ is connected,
    \item $G$ has at most one odd-degree vertex, and if it has none then it has at least one infinite-degree vertex,
    \item for every finite set of edges $F\subseteq E(G)$, $G - F$ has exactly one infinite connected component. 
\end{enumerate}
\end{theorem}

The characterization for directed graphs requires distinguishing whether an edge incident to a vertex $v$ is \textdef{incoming}, namely it is of the form $(u,v)$ for some vertex $u$, or \textdef{outgoing}, namely it is of the form $(v,u)$ for some vertex $u$. We write $\indeg_D(v)$ and $\outdeg_D(v)$ for the number of incoming and outgoing edges at $v$, respectively.

We say that a directed graph $D$ is \textdef{weakly unbiased} if, for every $X\subseteq V(D)$, the set of edges entering $X$ is infinite if and only if the set of edges leaving $X$ is infinite. We say that $D$ is \textdef{strongly unbiased}\footnote{In \cite{NashWilliams66}, the term \emph{unbiased} is used to refer to weak unbiasedness. Introducing the distinction between weak and strong unbiasedness makes some of the later proofs more readable.} if it is weakly unbiased and, whenever the sets of edges entering and leaving $X$ are both finite, their cardinalities differ by at most $1$. 

For each directed graph $D$, we define the \textdef{underlying undirected graph of $D$} as the undirected graph $H$ defined as follows: $V(H):=V(D) \cup \set{\xi_e\in \mathbb{N}}{e \in E(D)}$ where we can assume, modulo reindexing the vertices of $D$, that every $\xi_e$ is a fresh vertex not in $V(D)$. The edges of $H$ are obtained by splitting each edge $e\in E(D)$ in two using $\xi_e$ and then removing the orientation. In symbols: $E(H):= \bigcup_{e=(u,v)\in E(D)} \{ \{u, \xi_e \}, \{ \xi_e, v\} \}$. The reason why each edge is first split into two is to retain the fact that $D$ may contain the edges $(u,v)$ and $(v,u)$, hence removing the direction while preserving the edges would result in a multigraph. 

\begin{theorem}[Nash-Williams {\cite[Thm.\ 1]{NashWilliams66}}]
\thlabel{thm:nash-williams-one-way}
Let $D=(V,E)$ be a countable directed graph and let $v_0\in V$.  Then $D$ admits a one-way directed Eulerian path starting from $v_0$ iff all the following conditions hold:
\begin{enumerate}
	\item The underlying undirected graph of $D$ admits a one-way (undirected) Eulerian path,
	\item For every vertex $v\in V$,
	\[ \outdeg_D(v)=\begin{cases}
	\indeg_D(v)+1&\text{if }v=v_0,\\
	\indeg_D(v)&\text{if }v\neq v_0.
	\end{cases} \]
	with the convention that $\infty + 1 = \infty$.		
	\item $D$ is weakly unbiased. 
\end{enumerate}
\end{theorem}

We also mention the characterizations of when an undirected/directed graph admits a two-way Eulerian path.

\begin{theorem}[{Erd\H{o}s-Gr\"unwald-Weiszfeld \cite[Thm.\ 3.2.2]{Ore62}}]
	\thlabel{thm:EGW2}
A countable graph $G$ has an infinite two-way undirected Eulerian path iff 
\begin{enumerate}
    \item $G$ is connected,
    \item $G$ has no odd-degree vertex,
    \item for every finite set of edges $P$, $G- P$ has at most two infinite connected components. Moreover, if for any $v\in V(G)$, $P$ contains an even number of edges incident to $v$, then $G- P$ has exactly one infinite connected component.
\end{enumerate}
\end{theorem}

\begin{theorem}[Nash-Williams {\cite[Thm.\ 2]{NashWilliams66}}]
\thlabel{thm:nash-williams-two-way}
A countable directed graph $D=(V,E)$ admits a two-way directed Eulerian path iff all the following conditions hold:
\begin{enumerate}
	\item The underlying undirected graph of $D$ admits a two-way (undirected) Eulerian path,
	\item $D$ is strongly unbiased. 
\end{enumerate}
\end{theorem}

Observe that the two conditions in the previous theorem imply that, in a directed graph $D$ with a two-way directed Eulerian path, each vertex $v$ must satisfy $\indeg_D(v)=\outdeg_D(v)$.

\subsection{Weihrauch reducibility}\label{sec:weihrauch}
We briefly recall the main definitions on Weihrauch reducibility. For a more comprehensive presentation of Weihrauch reducibility, we refer the reader to \cite{BGP17}.

Given $f,g\pmfunction{\Baire}{\Baire}$, we say that $f$ is \textdef{Weihrauch reducible} to $g$, and write $f\weireducible g$, if there are two computable functionals $\Phi\pfunction{\Baire}{\Baire}$ and $\Psi\pfunction{\Baire\times\Baire}{\Baire}$ such that, for every $p\in\dom(f)$, 
\[ \Phi(p) \in \dom(g) \land (\forall q \in g(\Phi(p)))( \Psi(p,q ) \in f(p)). \] 

We call $\Phi$ and $\Psi$ the \textdef{forward} and \textdef{backward functionals}, respectively. These definitions can be extended to multi-valued functions on represented spaces by requiring that both $\Phi$ and $\Psi$ work on \emph{names} of the objects. However, spaces like $\mathbb{N}$, $\mathbb{Z}$, $\baire$, or $(\Baire)^\mathbb{N}$, as well as combinatorial objects like graphs or trees, admit a canonical representation as elements of $\Baire$. In these cases, with a small abuse of notation, we can identify an object with its name. In particular, as anticipated in the introduction, a graph $G=(V,E)$ is identified with the characteristic function of $V$ and $E$. One-way or two-way paths are canonically represented as elements of $\Baire$. We think of $\mathbb{N}\cup \{\infty\}$ as a discrete space (hence, in terms of representations, it can be identified with $\mathbb{N}$).

The Weihrauch degrees form a distributive lattice with a bottom element ($\emptyset$) and no top element. The join and meet operators can be obtained, respectively, by lifting the following two operators on multi-valued functions to the degree structure:
\begin{itemize}
	\item $(f_0\sqcup f_1)(i,x):=f_i(x)$, where $\dom(f_0\sqcup f_1):=\bigcup_{i<2} \{i\}\times \dom(f_i)$.
	\item $ (f_0\sqcap f_1)(x_0,x_1):= \bigcup_{i<2} \{i\}\times f_i(x_i)$, where $\dom(f_0\sqcap f_1):=\dom(f_0) \times \dom(f_1)$.
\end{itemize}
The framework of Weihrauch reducibility supports a wide range of operators on computational problems, enabling the study of the uniform computational strength of natural combinations and modifications of existing problems. We briefly introduce the ones that will be relevant for this paper. 

The \textdef{parallel product} $f\times g$ is the problem with $\dom(f\times g):= \dom(f)\times \dom(g)$ defined as $f\times g(x,y):= f(x)\times g(y)$. It corresponds to the idea of being able to use both $f$ and $g$ in parallel. The parallel product can be generalized to the countable case: we define the \textdef{parallelization} $\parallelization{f}$ of a problem $f$ as the problem whose domain is the set of sequences $\sequence{x_n}{n\in\mathbb{N}}$ of elements in $\dom(f)$ and $\parallelization{f}(\sequence{x_n}{n\in\mathbb{N}}):=\prod_{n\in\mathbb{N}} f(x_n)$. To capture the idea of using $f$ and $g$ in series, we introduce the \textdef{compositional product} $f * g$ as follows: let $\sequence{\Gamma_w}{w\in\Baire}$ be an effective presentation for all partial continuous functions with $\boldfacePi^0_2$ domain. We define  
\begin{gather*}
	\dom(f*g):=\set{(w,x)\in\Baire \times \dom(g)}{ (\forall y\in g(x))(\Gamma_w(y)\in \dom(f)) };\\
	f*g(w,x):=\set{(z,y)}{y\in g(x) \land z \in f(\Gamma_w(y))}.
\end{gather*}

The parallel product, parallelization, and compositional product lift to Weihrauch degrees. There are also many natural operators that are not degree-theoretic. The \textdef{jump}\footnote{This was originally introduced in \cite{BGM12}. Despite its name, $f'$ is not always strictly stronger than $f$ and does not lift to the Weihrauch degrees. See \cite{ALMMV25} for a more detailed discussion.} of a problem $f$ is the problem $f'$ whose input is a sequence $\sequence{p_n}{n\in\mathbb{N}}$ in $\Baire$ converging to some $p\in\dom(f)$, and $f'(\sequence{p_n}{n\in\mathbb{N}}):=f(p)$.  

For every partial multi-valued function $f\pmfunction{\Baire}{\Baire}$, its \textdef{total continuation} or \textdef{totalization} (\cite{NP18}) is the total multi-valued function $\totalization{f}\mfunction{\Baire}{\Baire}$
defined as
\[ \totalization{f}(x):=\begin{cases}
	f(x)	& \text{if } x\in\dom(f),\\
	\Baire	& \text{otherwise.}
\end{cases} \]
The \textdef{strong total continuation} of $f$ (\cite{MVRamsey}) is the total multi-valued function $\strongT{f}\mfunction{\Baire}{\Baire}$
defined as
\[ \strongT{f}(x):=\begin{cases}
	1\concat f(x)	& \text{if } x\in\dom(f),\\
	0\concat \Baire	& \text{otherwise.}
\end{cases} \]
The difference between $\totalization{f}$ and $\strongT{f}$ is that an answer for the latter comes with the indication of whether $x\in\dom(f)$, whereas this is not the case for the former. Neither of the two operators lifts to the Weihrauch degrees.

When studying the uniform computational strength of problems, it is common to compare them with \emph{benchmark} problems, whose Weihrauch degree is well-studied and can be used as a reference. 

For a metric space $X$, $\mflim_X$ is the problem that maps a convergent sequence in $X$ to its limit. We omit the subscript if $X=\Baire$. It is well-known that $\mflim_2 \strictlyweireducible \mflim_{\mathbb{N}} \strictlyweireducible \mflim$, where $2=\{0,1\}$ is the discrete space with two elements (\cite[Thm.\ 13.5]{BGM12}). A convenient way to use $\mflim$ is via the equivalence $\parallelization{\LPO}\weiequiv \mflim$ (\cite[Thm.\ 11.6.7]{BGP17}), where $\LPO\function{\Cantor}{2}$ is defined as $\LPO(p):=1$ iff $(\exists n)(p(n)=1)$. In other words, $\LPO$ can be used to answer a $\Sigma^0_1$ or $\Pi^0_1$ question, and $\mflim$ can solve countably many of these questions in parallel. Similarly, the problem $\LPO'$ can be used to answer a $\Sigma^0_2$ or $\Pi^0_2$ question, and $\mflim \strictlyweireducible \parallelization{\LPO'} \weiequiv \mflim'$ (see e.g.\ \cite[Thm.\ 11.6.5]{BGP17} and the proof of \cite[Thm.\ 7.2]{SV23}). 

The \emph{choice problems} are an important family of benchmarks: they intuitively describe the complexity of finding a point in a non-empty set. This complexity depends on the representation and the topological properties of the set. The formal definition of the \emph{closed choice problem} $\Choice{X}$ is usually given in the more general context of represented spaces, but for the purposes of this paper, we can give the following (Weihrauch-equivalent) definition: for $X=\{0,\hdots, k-1\}$ or $X=\mathbb{N}$ the problem $\Choice{X}$ takes as input a non-empty $A\subseteq X$ represented by a $\boldfacePi^0_1$ formula (i.e.\ the input is a relation $R\subseteq \mathbb{N}\times X$ such that $n\in A$ iff, for all $i\in\mathbb{N}$, $(i,n)\in R$), and produces $n\in A$. It is known that $\Choice{2}\strictlyweireducible \Choice{3} \strictlyweireducible \CNatural$ (see e.g.\ \cite[Prop.\ 11.7.18]{BGP17}). 

The problem $\CNatural$ is particularly important as being reducible to $\CNatural\weiequiv\mflim_\mathbb{N}$ (\cite[Thm.\ 11.7.13]{BGP17}) corresponds to being computable with finitely many mind-changes \cite[Prop.\ 11.6.1]{BGP17}. It is known that $\parallelization{\CNatural}\weiequiv \mflim$ \cite[Thm.\ 11.6.7]{BGP17}. In comparison, $\Choice{2}$ is significantly weaker, as $\parallelization{\Choice{2}}$ is Weihrauch-equivalent to \emph{weak K\"onig's lemma}, namely the problem of finding a path in an infinite tree $T\subseteq \cantor$. We will be especially interested in the problems $\Choice{2}'$ and $\CNatural'$. A convenient way to describe them is the following: for $X=2$ or $X=\mathbb{N}$, $\Choice{X}'$ is the problem that takes as input a non-empty $A\subseteq X$ defined by a $\boldfacePi^0_2$ formula, and produces $n\in A$. It is known that $\Choice{2}'$ is Weihrauch-equivalent to the pigeonhole principle $\mathsf{RT}^1_2$ for singletons and two colors (see \cite[Fact 2.3 and Prop.\ 3.4]{BRramsey17}).

Finally, we introduce the problem $\Sort_2\function{\Cantor}{\Cantor}$ defined as 
\[ \Sort_2(p):= \begin{cases}
0^n1^\omega & \text{if } p\text{ contains exactly $n$ zeroes} \\
0^\omega & \text{otherwise}. \\
\end{cases}\]
This will be useful to obtain an upper bound for the problem $\OrientOW$ (see \thref{def:orientow} below). It is known that $\CNatural \strictlyweireducible \Sort_2 \strictlyweireducible \mflim$ (\cite[Prop.\ 24]{NP18}).

\subsection{Some results}

We conclude this section with a few observations explaining why the techniques available in the literature do not directly yield a characterization of the uniform computational strength of finding an Eulerian path.

As mentioned in the introduction, every highly recursive or moderately computable graph that admits a (one-way or two-way) Eulerian path has a computable such path (\cite[Thms.\ 2 and 3]{Bean76} and \cite[Thm.\ A]{CarrascoVargas25}). In fact, both results relativize to $G$ and the degree function $d_G:=v \mapsto \deg_G(v)$. We also observe that Theorems 2 and 3 in \cite{Bean76} are uniform constructions, whereas the one-way case in \cite[Thm.\ A]{CarrascoVargas25} requires the additional knowledge of whether the graph has an odd-degree vertex. 

Let $\Deg$ be the function that takes as input a graph $G=(V,E)$ and a vertex $v\in V$ and returns $\deg_G(v)$. Let also $\DegLF$ be the restriction of $\Deg$ to locally finite graphs. 

\begin{proposition}
	\thlabel{thm:deg_locallyfinite}
	$\DegLF \weiequiv \CNatural \strictlyweireducible \Deg \weiequiv \sTCN$. 
\end{proposition}
\begin{proof}
	To show that $\DegLF \weireducible \CNatural$ and $\Deg \weireducible \sTCN$ we proceed as follows: for any graph $G=(V,E)$ and any vertex $v\in V$, first compute an enumeration $\sequence{v_n}{n\in\mathbb{N}}$ of $V$ (possibly with repetitions). Then, we can build the set 
	\[  A:=\set{b\in\mathbb{N}}{(\forall u \in V)( \{v,u\}\in E  \rightarrow (\exists n \le b)(u=v_n))}. \] 
	Observe that if $\deg_G(v)<\infty$ (which is always the case if $G$ is locally finite), then $A\neq \emptyset$, hence we can use $\CNatural$ to pick a point $b \in A$. Given such a bound $b$, we can compute $\deg_G(v)$ by computing the size of $\set{u \in V}{(\exists n\le b)(u=v_n)\land \{u,v\} \in E}$. Conversely, if $\deg_G(v)=\infty$ then $A=\emptyset$. In this case, the answer of $\sTCN$ contains the information that $A$ is empty, which implies that $\deg_G(v)=\infty$.

	To prove the reductions $\CNatural \weireducible \DegLF $ and $\sTCN\weireducible\Deg$ observe that, given $A\subseteq \mathbb{N}$, we can compute the function $g\function{\mathbb{N}}{\mathbb{N}}$ mapping $s$ to the least number not enumerated outside of $A$ by stage $s$. We can build a graph $G$ on $\mathbb{N}$ whose edges are all and only those of the form $\{0,s+1\}$ where $s$ is such that $g(s+1)\neq g(s)$. Observe that $A\neq\emptyset$ iff $g$ is eventually constant with value $\min A$ iff $\deg_G(0)$ is finite. Besides, $\deg_G(0)$ counts the changes of $g$, so, if this is a finite number, then it allows us to compute $\min A$.	

	Finally, it is known that $\CNatural\strictlyweireducible \sTCN$ as $\parallelization{\CNatural}\weiequiv\mflim \strictlyweireducible  \mflim'\weiequiv \parallelization{\sTCN}$ 
	(the last equivalence follows e.g.\ from \cite[Thm.\ 7.12]{SV23}).
\end{proof}

The above proposition can be adapted so that, in the proof of the reductions $\CNatural \weireducible \DegLF $ and $\sTCN\weireducible\Deg$, we always build a graph $G$ that admits a one-way or two-way Eulerian path. 

\begin{corollary}
	\thlabel{thm:EP_bounded}
	The restriction of $\EP$ to locally finite graphs is Weihrauch-equivalent to $\mflim$.
\end{corollary}
\begin{proof}
	For the upper bound, let $G=(V,E)$ be a locally finite graph admitting a one-way Eulerian path. In light of \cite[Thm.\ 2]{Bean76}, in order to compute an Eulerian path through a given locally finite graph, it is enough to compute the function $d\function{V}{\mathbb{N}}$ mapping every vertex to its degree. By the previous proposition, this can be done uniformly using $\mflim$.

	The converse reduction is implicit in the direction $(3)\Rightarrow (1)$ of the proof of \cite[Thm.\ 12]{GH98}.
\end{proof}

With minor modifications, the above arguments can be adapted to the case of directed graphs. In fact, \thref{thm:deg_locallyfinite} can be adapted to show that the map $(D,v)\mapsto \indeg_D(v)$ is Weihrauch-equivalent to $\CNatural$ if each vertex has a finite number of incoming edges, and it is Weihrauch-equivalent to $\sTCN$ in the general case. An analogous characterization holds for the map $(D,v)\mapsto \outdeg_D(v)$. In particular, knowing that for every $D\in\dom(\DEP)$ and every $v\in V(D)$, $\indeg_D(v)$ is finite iff $\outdeg_D(v)$ is finite, we obtain that the restriction of $\DEP$ to locally finite graphs is Weihrauch-equivalent to $\mflim$.

These observations show that any approach to bounding the uniform computational strength of $\EP$ or $\DEP$ that begins by computing the degree function (or both the in-degree and out-degree functions in case of directed graphs) would yield an upper bound at the level of $\mflim'$ (or above). In light of the rest of the paper, such an upper bound would not be optimal.

Similar considerations apply to the case of two-way undirected Eulerian paths. Indeed, as noted in \cite[Sec.\ 3]{GH98}, the argument for the one-way Eulerian paths can be adapted to two-way undirected Eulerian paths. An analogous argument for two-way directed Eulerian paths does not seem to have been explicitly presented in the literature.

\section{One-way Eulerian paths}

For the sake of readability, we split the characterization of the Weihrauch degree of $\EP$ and $\DEP$ as follows: we first show the lower bound (\thref{lem:lim-choice-to-ep}). Then, for the upper bound, we show that starting from a prescribed vertex simplifies the problem of finding a one-way (directed or undirected) Eulerian path (\thref{thm:fixed-start}). These facts, together with a technical lemma (\thref{lem:lim-first-order}) are then combined in the proof of \thref{thm:EP}.

\begin{lemma}
	\thlabel{lem:lim-choice-to-ep}
    $\mflim\times\CNatural'\weireducible\EP$ and $\mflim\times\CNatural'\weireducible\DEP$.
\end{lemma}
\begin{proof}
	We only show that $\mflim\times\CNatural'\weireducible\DEP$. The reduction to $\EP$ builds the same graph, with the only difference being that the edges are undirected. 

	Let $\sequence{p_n}{n\in\mathbb{N}}$ be an input for $\parallelization{\LPO}$ ($\weiequiv\mflim)$ and $A\neq \emptyset$ be a $\boldfacePi^0_2$ subset of $\mathbb{N}$. We represent $A$ as a sequence $\sequence{q_n}{n\in\mathbb{N}}$ in $\Cantor$ such that $n\in A \iff (\exists^\infty i)(q_n(i)=1)$. 

	We define a graph $D=(V,E)$ as follows (see Figure~\ref{fig:2gadget}). We use disjoint computable sets of labels for the backbone, the vertices $h_n$, and auxiliary vertices. The graph has a ray-like backbone $b_0, b_1, \hdots$ where every backbone vertex $b_i$ is connected to $b_{i+1}$ via two paths of length 2 going in opposite directions.
	
	We will attach two different types of gadgets to this backbone. One will help us decide each $\LPO$ instance and the other will give us a solution to $\CNatural'(A)$.

	The first one is (essentially) the same gadget used in \cite[Thm.\ 12]{GH98}: when we first see that $p_e(s)=1$ for some $s$, we pick two fresh indexes $u_e$ and $v_e$ and attach a triangle to $b_{2e+1}$, oriented as $b_{2e+1}\to u_e \to v_e \to b_{2e+1}$.

	To solve the $\CNatural'$-instance, for each $n$ we define a dedicated vertex $h_n$ (we choose $h_n$ so that, knowing $h_n$, we can uniformly compute $n$) and attach it to a sufficiently large even backbone vertex unused so far by a triangle $h_n \to c_{n,j}\to b_{2k} \to h_n$, where each $c_{n,j}$ is new. Whenever we find a new $i$ such that $q_n(i)=1$, we add another such triangle, using a fresh auxiliary vertex $c_{n,j}$ and a sufficiently large backbone vertex $b_{2k}$ not used by any $h_m$ so far. 
	Notice that $\indeg_D(h_n)=\outdeg_D(h_n)=\infty$ iff $n \in A$. 
	\begin{figure}[htbp]
	\centering
	\begin{tikzpicture}[x=1.05cm,y=0.85cm,>=stealth,
	vertex/.style={circle,fill=black,inner sep=1.6pt}]
	\foreach \i in {0,1,2,4,5,6} {
	\node[vertex] (b\i) at (\i,0) {};
	\node[below=3pt] at (b\i) {$b_{\i}$};
	}
	\node[vertex] (b3) at (3,0) {};
	\node[above=3pt] at (b3) {$b_3$};
	\foreach \i in {0,...,5} {
	\pgfmathtruncatemacro{\j}{\i+1}
	\node[vertex] (x\i) at (\i+0.5,0.28) {};
	\node[vertex] (y\i) at (\i+0.5,-0.28) {};
	\draw[->] (b\i)--(x\i);
	\draw[->] (x\i)--(b\j);
	\draw[->] (b\j)--(y\i);
	\draw[->] (y\i)--(b\i);
	}
	\node at (6.5,0) {$\cdots$};
	\node[vertex,label=above:{$h_n$}] (h) at (2,2.5) {};
	\foreach \i in {0,2,4} {
	\node[vertex] (c\i) at (\i-0.3,1.1) {};
	\draw[->] (h)--(c\i);
	\draw[->] (c\i)--(b\i);
	\draw[->] (b\i)--(h);
	}
	\node[left=3pt] at (c0) {$c_{n,0}$};
	\node[left=0.1pt] at (c2) {$c_{n,1}$};
	\node[right=3pt] at (c4) {$c_{n,2}$};
	\node at (5,1.1) {$\cdots$};
	\node[vertex,label=left:{$u_1$}] (u) at (2.5,-1.4) {};
	\node[vertex,label=right:{$v_1$}] (v) at (3.5,-1.4) {};
	\draw[->] (b3)--(u);
	\draw[->] (u)--(v);
	\draw[->] (v)--(b3);
	\end{tikzpicture}
	\caption{The construction of the directed graph $D$.}
		\label{fig:2gadget}
	\end{figure}

	Let us now show that $D$ admits a one-way directed Eulerian path. 
	Let $G$ be the underlying undirected graph of $D$. Observe that $\deg_G(v)=\infty$ iff $v\in \set{h_n}{n\in A}$, and otherwise it has even degree. In fact, for every $v\in V(D)$, $\indeg_D(v)=\outdeg_D(v)$. Besides, the fact that $A\neq\emptyset$ guarantees that at least one infinite-degree vertex exists. Observe also that deleting any finite set $F\subseteq E(G)$ of edges leaves $G$ with exactly one infinite component: indeed, for a sufficiently large backbone vertex $b$, the backbone tail beyond $b$ remains connected, and each backbone vertex before $b$ can have at most $1$ triangle attached (whether a connection to some $h_n$ or a witness for the $e$-th $\LPO$ instance). By \thref{thm:EGW}, this shows that $G$ admits a one-way Eulerian path. It remains to show that $D$ is weakly unbiased. Observe first that its edges can be partitioned into finite directed cycles (either triangles or diamonds). For every $X\subseteq V(D)$, each such cycle has equally many edges entering and leaving $X$. Summing over the cycles shows that the sets of edges entering and leaving $X$ have the same cardinality, hence $D$ is weakly unbiased. In particular, for every $n\in A$, \thref{thm:nash-williams-one-way} implies that there is a one-way directed Eulerian path through $D$ starting from $h_n$.

	Let $P\in \DEP(D)$ be a one-way Eulerian path through $D$. Since every vertex has the same incoming and outgoing degree, $P(0)=h_n$ for some $n\in A$. This is because the starting vertex of the path must satisfy $\indeg_D(P(0))+1=\outdeg_D(P(0))$. Therefore knowing $P(0)$ lets us uniformly compute an answer for $\CNatural'(A)$. 

	To decide whether $\LPO(p_e)=1$, we wait until $P$ has traversed all four backbone edges incident to $b_{2e+1}$. Then $\LPO(p_e)=1$ iff the path traverses the triangle $b_{2e+1} \to  u_e \to  v_e \to b_{2e+1}$. %This shows that we can uniformly compute the answer to every $\LPO$-instance, and therefore conclude the proof.
	Indeed, the existence of such triangle implies that we have found $s$ such that $p_e(s)=1$. Conversely, if the path does not traverse the triangle, it means it has never been added: if the triangle is added to the graph but not traversed after visiting the four edges of $b_{2e+1}$ then it would be disconnected from the rest of the graph, hence unreachable, against the fact that $P$ is an Eulerian path. This shows that we can uniformly compute the answer to every $\LPO$-instance, and therefore conclude the proof.
\end{proof}

Observe that, in the previous argument, the answer to the $\CNatural'$-instance is determined solely by the initial vertex of the Eulerian path. This is no coincidence: the next result shows that specifying a valid starting vertex leaves a problem solvable using $\mflim$ alone.

\begin{theorem}
\thlabel{thm:fixed-start}
	Let $\EPstart$ (resp.\ $\DEPstart$) be the problem that takes as input a pair $(G,v)$, where $G=(V,E)$ is an undirected (resp.\ directed) graph, $v\in V$, and $G$ admits a one-way Eulerian path starting from $v$, and returns a one-way undirected (resp.\ directed) Eulerian path through $G$ starting from $v$. Then
	\[ \EPstart\weiequiv  \DEPstart \weiequiv \mflim. \]
\end{theorem}
\begin{proof}
	The reductions $\mflim\weireducible \EPstart$ and $\mflim\weireducible \DEPstart$ can be obtained by simplifying the construction in the proof of \thref{lem:lim-choice-to-ep} as follows (see Figure~\ref{fig:directed-fixed-start-ray}): given an instance $\sequence{p_n}{n\in\mathbb{N}}$ for $\parallelization{\LPO}$, define a directed graph $D$ so that it contains a ray $b_0 \to b_1 \to b_2 \hdots$ (we connect $b_i$ and $b_{i+1}$ directly instead of via a diamond) and, when we first see that $p_e$ contains a $1$, we add a triangle $b_e \to u_e \to v_e \to b_e$, for some nodes $u_e,v_e$ with sufficiently large index. The construction is the same in the undirected case: simply remove the orientation of the edges. 

	\begin{figure}[htbp]
\centering
\begin{tikzpicture}[x=1.05cm,y=0.85cm,>=stealth,
 vertex/.style={circle,fill=black,inner sep=1.6pt}]
 \foreach \i in {0,...,6} {
  \node[vertex,label=above:{$b_{\i}$}] (b\i) at (\i,0) {};
 }
 \foreach \i in {0,...,5} {
  \pgfmathtruncatemacro{\j}{\i+1}
  \draw[->] (b\i)--(b\j);
 }
 \draw[->] (b6)--(6.4,0);
 \node at (6.65,0) {$\cdots$};
 \foreach \i in {1,4} {
  \node[vertex,label=left:{$u_{\i}$}] (u\i) at (\i-0.45,-1.1) {};
  \node[vertex,label=right:{$v_{\i}$}] (v\i) at (\i+0.45,-1.1) {};
  \draw[->] (b\i)--(u\i);
  \draw[->] (u\i)--(v\i);
  \draw[->] (v\i)--(b\i);
 }
 \node[left=5pt] at (b0) {start};
\end{tikzpicture}
\caption{The construction for the reduction $\mflim\weireducible\DEPstart$.}
\label{fig:directed-fixed-start-ray}
\end{figure}
	It is apparent that $D$ admits a one-way Eulerian path starting from $b_0$. Besides, we can find the solution to $\LPO(p_e)$ by checking if the path $P\in \DEPstart(D,b_0)$ visits $b_{e+1}$ immediately after $b_e$ or if it instead visits the triangle $b_e \to u_e \to v_e \to b_e$ first. 

	Let us now show that $\EPstart\weireducible \mflim$. Let $(G,v_0)\in\dom(\EPstart)$. As a consequence of \thref{thm:EGW}, a finite trail $\sigma$ starting at $v_0$ is extendible to a one-way Eulerian path iff the final vertex of $\sigma$ has an unused incident edge and $G-E(\sigma)$ is connected (see \cite[Thm.\ 3]{CarrascoVargas25}). We show that these two conditions are uniformly $\Delta^0_2$ relative to $\sigma$ and $G$, hence we can use $\mflim$ to check them in parallel for every $\sigma$. Fix a finite trail $\sigma$ and assume without loss of generality that $\length{\sigma}\ge 2$. Checking if $w:=\sigma(\length{\sigma}-1)$ has an unused incident edge in $G-E(\sigma)$ is clearly $\Sigma^0_1$ relative to $\sigma$ and $E$. For the second condition, let
	\[ B_\sigma := \set{b \in \ran(\sigma)}{(\exists x \in V(G))(\{x,b\}\in E(G)\setminus E(\sigma))}  \]
	be the set of vertices in $\sigma$ that are the endpoints of some edge not in $E(\sigma)$.
	Since $G$ is connected, every connected component of $G-E(\sigma)$ contains a vertex in $B_\sigma$. In particular, $G-E(\sigma)$ is connected iff $B_\sigma$ is connected in $G-E(\sigma)$. Since $B_\sigma$ is $\Sigma^0_1$, the latter condition is $\Delta^0_2$ relative to $\sigma$ and $G$.

	Using $\mflim$, we can then uniformly compute a one-way Eulerian path through $G$ inductively by unbounded search: fix an enumeration $\sequence{e_n}{n\in\mathbb{N}}$ of $E(G)$ and let $\sigma_0:=\str{v_0}$. At each stage $s$, we search for some finite trail $\sigma_{s+1} \pextend \sigma_s$ that is extendible to an Eulerian path and that visits $e_s$. By the inductive hypothesis, such $\sigma_{s+1}$ exists and the information given to us by $\mflim$ allows us to find one. The limit of $\sequence{\sigma_s}{s\in\mathbb{N}}$ gives us a one-way Eulerian path through $G$.

	The proof of $\DEPstart\weireducible \mflim$ follows a similar strategy: let $(D,v_0)$ be an instance of $\DEPstart$. Using Nash-Williams's theorem (\thref{thm:nash-williams-one-way}), we can show that a finite directed trail $\sigma$ starting from $v_0$ is extendible iff 
	\begin{itemize}
		\item the final vertex $w:=\sigma(\length{\sigma}-1)$ has an unused outgoing edge and
		\item the underlying undirected graph of $H:=D-E(\sigma)$ is connected. 
	\end{itemize}
	Clearly, both conditions are necessary. To see that they are sufficient, observe that they imply that the underlying undirected graph of $H$ has a one-way Eulerian path (as above) and removing finitely many edges preserves weak unbiasedness. Since $\sigma$ is a directed trail, for every vertex $v\neq v_0, w$ there is $k$ such that $\outdeg_D(v)=\outdeg_H(v)+k$ and $\indeg_D(v)=\indeg_H(v)+k$. The same holds for $v_0$ if $v_0=w$ (hence, in particular, we obtain $\outdeg_H(w) = \indeg_H(w) +1$). Otherwise, $\outdeg_H(v_0) = \indeg_H(v_0)$ (as the first edge goes out of $v_0$ but there is no previous edge going to $v_0$) and 
	\[ \outdeg_D(w)=\outdeg_H(w)+k \quad \text{ and } \quad \indeg_D(w)=\indeg_H(w)+k +1. \]
	This shows that $H$ and $w$ satisfy the hypotheses of \thref{thm:nash-williams-one-way} and hence $H$ admits a directed one-way Eulerian path starting at $w$.

	The rest of the proof follows as in the undirected case: the two conditions above are uniformly $\Delta^0_2$ relative to $\sigma$ and $D$, hence we can use $\mflim$ to check them in parallel for every $\sigma$. Then, we can uniformly compute a directed one-way Eulerian path through $D$ starting from $v_0$ proceeding inductively by unbounded search while ensuring that we visit each edge.
\end{proof}

\begin{lemma}
	\thlabel{lem:lim-first-order}
    For every problem $g$ with codomain $\mathbb{N}$, $\mflim * g\weiequiv\mflim\times g$.
\end{lemma}
\begin{proof}
Without loss of generality, assume $g\pmfunction{\Baire}{\mathbb{N}}$. The reduction $\mflim\times g \weireducible \mflim * g$ is immediate (\cite[Prop.\ 11.5.7]{BGP17}), hence we only need to show the converse reduction. 
Let $(w,x)\in \Baire\times \Baire$ be a valid input for $\mflim * g$. We show that we can solve it using $\parallelization{\LPO}\times g$. We use the $\pairing{i,j,s}$-th instance of $\LPO$ to check if there is $t>s$ such that, in $t$ steps, $\Gamma_w(i)$ produces a prefix of a sequence $\sequence{p_n}{n\in\mathbb{N}}$ such that for some $s<m<t$, $p_s(j)\neq p_m(j)$. 

Notice $\Gamma_w(i)$ need not produce a convergent sequence if $i\notin g(x)$. On the other hand, if $i\in g(x)$ then $\Gamma_w(i)=\sequence{p_n}{n\in\mathbb{N}}\in\dom(\mflim)$ and, for every $j$ there is $s$ such that $\sequence{p_n(j)}{n\in\mathbb{N}}$ stabilizes after stage $s$. In particular, knowing some solution $i$ for $g(x)$, we can uniformly compute $\mflim(\Gamma_w(i))$ by searching, for each $j$, for the first $s$ such that the $\pairing{i,j,s}$-th instance of $\LPO$ returns $0$. 
\end{proof}

We are now ready to characterize the Weihrauch-degree of $\EP$ and $\DEP$.

\begin{theorem}
\thlabel{thm:EP}
$\EP\weiequiv\DEP\weiequiv\mflim\times\CNatural'$.
\end{theorem}
\begin{proof}
The fact that $\mflim\times\CNatural'$ is a lower bound for both problems is given by \thref{lem:lim-choice-to-ep}. The two reductions $\EP\weireducible\mflim\times\CNatural'$ and $\DEP\weireducible\mflim\times\CNatural'$ follow a similar strategy: we first use $\CNatural'\weiequiv \CNatural' * \LPO'$ (as it is known that $\LPO'\weireducible \CNatural'$ and $\CNatural'$ is closed under compositional product \cite[Thm.\ 11.7.53(1)]{BGP17}) to find a starting vertex for a (one-way) Eulerian path. We then use $\mflim$ to obtain the rest of the path (using \thref{thm:fixed-start}). This shows that both $\EP$ and $\DEP$ reduce to $\mflim*\CNatural'$, and the claim follows by \thref{lem:lim-first-order}.

In the undirected case, let $G\in\dom(\EP)$. Observe that checking if $\deg_G(v)$ is odd is a uniformly $\Sigma^0_2$ condition relative to $G$ and $v$, as it is equivalent to 
	\[(\exists n\in\mathbb{N}) (\exists F\subseteq E, |F|=2n+1) (\forall e\in E)(e \text{ is incident to }v \leftrightarrow e \in F)\] 
If there is a vertex with odd degree, then use $\CNatural'$ to choose the unique such vertex (as $\CNatural$ suffices to choose an element from a given $\boldfaceSigma^0_2$ subset of $\mathbb{N}$). Otherwise, there must be a vertex with infinite degree, hence we use $\CNatural'$ to choose from the set of vertices with infinite degree (which is uniformly $\Pi^0_2$ relative to $G$). 
	
In the directed case, let $D\in\dom(\DEP)$. We can use $\LPO'$ to check if there is a vertex $v\in V$ such that $\outdeg_D(v)=\indeg_D(v)+1<\infty$. Indeed, this condition is equivalent to check if there are a vertex $v$ and two finite sets of edges $F_{\mathrm{in}},F_\mathrm{out}\subset E$ such that $|F_{\mathrm{out}}| = |F_{\mathrm{in}}| +1$, and $F_{\mathrm{in}}$ (resp.\ $F_{\mathrm{out}}$) contains all and only the edges in $E(D)$ incoming to (resp.\ outgoing from) $v$. In particular, this is uniformly $\Sigma^0_2$ relative to $D$. If such a vertex exists, we can use $\CNatural'$ to find it. Otherwise, we can use $\CNatural'$ to choose from the set of vertices $v\in V(D)$ such that $\outdeg_D(v)=\indeg_D(v)=\infty$ (which is uniformly $\Pi^0_2$ relative to $D$).
\end{proof}

\begin{corollary}
	$\mflim'\not\weireducible \EP$ and $\mflim'\not\weireducible \DEP$.
\end{corollary}
\begin{proof}
	Simply follows from the fact that there is a computable instance of $\mflim'$ whose solution is $\emptyset''$, while every computable instance of $\mflim \times \CNatural'\weiequiv \EP$ has a $\Delta^0_2$ solution.
\end{proof}

\begin{corollary}
	Both $\EP$ and $\DEP$ are closed under parallel product, but not parallelizable nor closed under compositional product.
\end{corollary}
\begin{proof}
	Since both $\mflim$ and $\CNatural'$ are closed under parallel product, so are $\EP$ and $\DEP$. To show that $\EP$ and $\DEP$ are not closed under parallelization nor compositional product, it is enough to observe that $\parallelization{\CNatural'} \weiequiv \mflim' \weiequiv \mflim * \mflim$, which can be proved using e.g.\ \cite[Props.\ 11.6.12 and 11.6.14]{BGP17}.
\end{proof}

\section{Orienting the edges}
\label{sec:orientations}

\begin{definition}
	\thlabel{def:orientow}
	Let $\OrientOW$ be the problem that, given $G=(V,E)\in\dom(\EP)$, returns $D=(V,A)\in\dom(\DEP)$ such that the map $\pi\function{A}{E}:=(v,w) \mapsto \{v, w\}$ is a (computable) bijection.
\end{definition}

It is not hard to show that $\OrientOW$ is not uniformly computable.

\begin{proposition}
	\thlabel{thm:lim2<=orient}
	$\mflim_2\weireducible \OrientOW$. 
\end{proposition}
\begin{proof}
	Let $\sequence{x_n}{n\in\mathbb{N}}$ be an instance of $\mflim_2$. We build a graph $G=(\mathbb{N},E)$ as follows: at stage $0$, let $l_0:=0$, $r_0:=1$ and we add the edge $\{0,1\}$ to $E$. Intuitively, we are building a ray, where the direction of the ray (left or right) is determined by the current element of the input sequence. At stage $n+1$, if $x_{n+1}=0$, we add the edge $\{n+2, l_n\}$ and define $l_{n+1}:=n+2, r_{n+1}:=r_n$. Symmetrically, if $x_{n+1}=1$, we add the edge $\{n+2, r_n\}$ and define $l_{n+1}:=l_n, r_{n+1}:=n+2$. 
	
	Since the sequence $\sequence{x_n}{n\in\mathbb{N}}$ is eventually constant with value $i$, exactly one side is extended infinitely often, so $G$ is a ray. The only $D\in\OrientOW(G)$ orients $\{0,1\}$ as $(0,1)$ iff $i=1$, uniformly computing the limit.
\end{proof}

We can obtain the separation $\OrientOW\not\weireducible \mflim_2$ as a corollary of the next stronger result (as $\mflim_2\strictlyweireducible \CNatural$): 

\begin{proposition}
	\thlabel{thm:orient_CN}
	$\OrientOW\not\weireducible \CNatural$.
\end{proposition}
\begin{proof}
	Assume towards a contradiction that $\OrientOW\weireducible \CNatural$ via the functionals $\Phi,\Psi$. Notice that, for each input $A$ for $\CNatural$, we can uniformly compute a guessing function $g_A\function{\mathbb{N}}{\mathbb{N}}$ such that, for each $s$, $g_A(s)$ is the least number not enumerated outside of $A$ by stage $s$. In particular, $\lim_s g_A(s) = \min A$. 

	To obtain a contradiction, we build a graph $G$ so as to force the guessing function $g_{\Phi(G)}$ to not have a limit.  The graph $G$ is built iteratively as follows: at each stage $s$, the approximation $G_s$ of $G$ is a graph consisting of a finite number of finite cycles going through the vertex $0$ and a finite ray $R_s$ containing $0$. Let $g_s$ be the least number not enumerated outside of $\Phi(G_s)$ by stage $s$. Without loss of generality, we can assume that, in $s$ steps, both computations $\Phi(G_s)$ and $\Psi(G_s,g_s)$ only query the oracle on edges whose endpoints are at most $s$. This is useful so that later changes to the approximation of $G$ preserve earlier computations.
	
	We start at stage $1$ with $G_1$ being the graph with two vertices $0$ and $1$ connected by an edge. The construction alternates between two different phases. We start in Phase 1. 
	\begin{description}
		\item[Phase 1] In this phase, the ray $R_s$ is of the form $\str{0,t,t+1,t+2,\hdots,s}$ for some $t\le s$. We check if $\Psi(G_s,g_s)$, in $s$ steps, commits to the orientation of any edge $\{ t+j,t+j+1 \}$ as $(t+j,t+j+1)$. If this does not happen, we add the edge $\{s+1,s\}$, go to the next stage and remain in Phase 1. Otherwise, we add the edge $\{s+1,0\}$ and we move to Phase 2 and go to the next stage.
		\item[Phase 2] The ray is $R_s=\str{s,\ldots,r,0,t,\ldots,t+k}$ for some $t<t+k<r\le s$. We add the edge $\{s+1,s\}$. We also check if $\Psi(G_s,g_s)$, in $s$ steps, commits to the orientation of any edge $\{r,r+1\}$ as $(r,r+1)$. If this does not happen, we go to the next stage and remain in Phase 2. Otherwise, we close the loop by adding the edge $\{s+1, t+k\}$. We also add the edge $\{0,s+2\}$, then move to Phase 1 and go to stage $s+2$.
	\end{description}
		\begin{figure}[htbp]
	\centering
	\begin{tikzpicture}[x=0.78cm,y=0.78cm,
		vertex/.style={circle,fill=black,inner sep=1.4pt},
		every node/.style={font=\small}]
		% Phase 1: extend the right arm.
		\begin{scope}
			\node at (0,1.25) {Phase 1};
			\node[vertex,label=below:{$0$}] (a0) at (-1.5,0) {};
			\node[vertex] (a1) at (-0.6,0) {};
			\node[vertex] (a2) at (0.3,0) {};
			\node[vertex,label=below:{$s$}] (a3) at (1.2,0) {};
			\draw (a0)--(a1)--(a2)--(a3);
			\draw[dashed] (a3)--(1.8,0);
			\draw[->] (-0.5,0.55)--(1.2,0.55);
			\node at (0,-1.15) {Extend to the right};
		\end{scope}
		% Phase 2: keep the right endpoint fixed and extend left.
		\begin{scope}[xshift=4.5cm]
			\node at (0,1.25) {Phase 2};
			\node[vertex,label=below:{$0$}] (b0) at (0,0) {};
			\node[vertex] (b1) at (0.8,0) {};
			\node[vertex,label=below:{$t+k$}] (b2) at (1.6,0) {};
			\node[vertex] (b3) at (-0.8,0) {};
			\node[vertex,label=below:{$s$}] (b4) at (-1.6,0) {};
			\draw (b0)--(b1)--(b2) (b0)--(b3)--(b4);
			\draw[dashed] (b4)--(-2,0);
			\draw[->] (b1)--(b2);
			\draw[->] (-0.2,0.55)--(-1.6,0.55);
			\node at (0,-1.15) {Extend to the left};
		\end{scope}
		% Close the two arms into a cycle.
		\begin{scope}[xshift=9cm]
			\node at (0,1.25) {Close the loop};
			\node[vertex,label=below:{$0$}] (c0) at (0,0) {};
			\node[vertex] (c1) at (0.8,0) {};
			\node[vertex,label=below:{$t+k$}] (c2) at (1.6,0) {};
			\node[vertex] (c3) at (-0.8,0) {};
			\node[vertex,label=below:{$s+1$}] (c4) at (-1.6,0) {};
			\draw (c0)--(c1)--(c2) (c0)--(c3)--(c4);
			\draw[->] (c1)--(c2);
			\draw[->] (c3)--(c4);
			\draw[dashed] (c4) to[out=65,in=115] (c2);
			\node at (0,-1.15) {Join the endpoints};
		\end{scope}
	\end{tikzpicture}
	\caption{The two phases of the construction and the closing of the loop. Previously completed loops at $0$ are omitted. Arrows on edges indicate the orientations produced by $\Psi$.}
	\label{fig:orient-diagonalisation}
	\end{figure}

	Observe that the construction must alternate between the two phases infinitely often. Indeed, if from a certain stage $s$ on, the construction remains in Phase $1$, then the resulting graph $G$ is a ray with finitely many attached cycles (hence $G\in\dom(\EP)$) and the backward functional never orients the edges of the ray away from the origin $0$. Clearly, no edge of the ray can be oriented toward $0$ (otherwise there cannot be any Eulerian path), hence it means that $\Psi$ never produces an orientation for any edge of the ray, against the fact that $\Psi$ witnesses a Weihrauch reduction. The case when the construction remains in Phase 2 is symmetric. 

	This shows that the resulting graph $G$ is the union of countably many cycles whose only vertex in common is $0$ (which therefore has infinite degree), therefore $G$ admits an Eulerian path.

	To conclude the proof, it is enough to notice that every change from Phase 1 to Phase 2 forces $g_{\Phi(G)}$ to change guess. Indeed, assume otherwise and let $g_s$ be the current guess. When we close the loop, $\Psi(G,g_s)$ has committed to orienting an edge away from $0$ on each of its two sides. Therefore, regardless of how the remaining edges are oriented, there must be a vertex in the loop other than $0$ with two incoming edges and no outgoing one. Since this vertex has degree $2$ in $G$, the resulting directed graph cannot admit an Eulerian path, against the fact that $\Psi$ witnesses a Weihrauch reduction. Hence the guesses must change infinitely often, contradicting the convergence of $g_s$.
\end{proof}

We now work towards producing an upper bound for $\OrientOW$. We will show that $\OrientOW \strictlyweireducible \Choice{2}' \sqcap \Sort_2$. To show that $\OrientOW \weireducible \Choice{2}'$, we first need the following lemma.

\begin{lemma}
	\thlabel{thm:D_Dop}
	There is a computable procedure that, given a graph $G=(V,E)\in \dom(\EP)$, computes a directed graph $D=(V,A)$ by choosing an orientation of the edges in $E$ such that either $D\in\dom(\DEP)$ or $D^\mathrm{op}:=(V,A^\mathrm{op})\in \dom(\DEP)$.
\end{lemma}
\begin{proof}
	The idea is to orient the graph locally, partitioning its edges into finite directed closed trails and at most one remaining directed trail (which could be finite, one-way infinite, or two-way infinite). The hypothesis that $G$ admits a one-way Eulerian path guarantees that the searches in the construction succeed. We will then show that the resulting orientation is weakly unbiased and that every vertex is balanced, with the possible exception of the unique odd-degree vertex of $G$. If there is such a vertex $v$, the correct global direction is the one satisfying $\indeg_D(v)+1=\outdeg_D(v)$; otherwise, both $D$ and $D^{\mathrm{op}}$ will admit a one-way directed Eulerian path.

	We now show how to define the orientation $A$ (which yields the directed graph $D$). Let $\sequence{e_n}{n\in\mathbb{N}}$ be an enumeration of $E$. At each step, we maintain a finite directed open trail $\tau_s$, which we aim to extend at either endpoint. Start with $\tau_0$ being an arbitrary orientation of $e_0$. At stage $s$, let $e_n$ be the least unoriented edge yet. We search for a finite directed trail $\rho$ that only includes unoriented edges and that satisfies one of the following:
\begin{enumerate}
	\item $\rho$ includes $e_n$ and is a closed trail edge-disjoint from $\tau_s$;
	\item $\rho$ includes $e_n$ and can be prepended or appended to $\tau_s$ to form an open directed trail;
	\item $\rho$ can be appended to $\tau_s$ to form a closed trail.
\end{enumerate}
In case (1), $\rho$ defines an orientation of its edges. We define $\tau_{s+1}:=\tau_s$ and go to the next stage. Similarly, in case (2), we define $\vartheta$ as the directed trail obtained by prepending or appending $\rho$ to $\tau_s$. Then, we define $\tau_{s+1}:=\vartheta$ and go to the next stage. %Otherwise, as in (1), we use $\rho$ to define an orientation of $\vartheta$. We define $\tau_{s+1}:=\str{x_m, y_m}$, where $e_m = \{ x_m, y_m \}$ is the first edge not oriented so far and go to the next stage. 
Finally, in case (3), we have obtained a closed directed trail. As in (1), we orient the mentioned edges as specified by the obtained closed trail. Then, we define $\tau_{s+1}:=\str{x_m, y_m}$, where $e_m = \{ x_m, y_m \}$ is the first edge not oriented so far, and go to the next stage.

We now show that the search succeeds. Suppose that, at stage $s$, none of (1),(2), or (3) is possible. Let $u:=\tau_s(0)$ and $w:=\tau_s(\length{\tau_s}-1)$. Let $H:=(V(G),E_s)$, where $E_s$ is the set of edges not oriented by stage $s$. Observe that committing to the orientation of a finite closed trail preserves the degree parity of each vertex. Similarly, the degree parity is preserved for every vertex in $\tau_s$ different from $u,w$. In symbols:
\[ (\forall v\in V\setminus\{u,w\})(\exists k\in\mathbb{N})( \deg_G(v)=\deg_H(v)+2k).\] 
Let $K$ be the connected component of $H$ containing $e_n$. 

The failure of (2) and (3) implies that $K$ contains neither $u$ nor $w$. Besides, the failure of (1) implies that $K$ is infinite. Indeed, recall that $G$ can have at most one odd-degree vertex. If $K$ were finite then all of its vertices would have even degree, since the sum of the degrees of vertices in a finite graph is always even\footnote{This is an elementary fact, see e.g.\ \cite[Eq.\ (1.2.5), p.\ 7]{Ore62}}. In particular, $K$ would admit a closed Eulerian trail \cite[Thm.\ 3.11]{Ore62}. Such a trail, once oriented, would contradict the failure of (1). 

By \thref{thm:EGW}, $K$ is the unique infinite connected component of $H$. Let $K_u$ and $K_w$ be the connected components in $H$ of $u$ and $w$ respectively. In particular $K_u \neq K \neq K_w$.  Assume towards a contradiction that $K_u\neq K_w$. Within $K_u$, deleting the already oriented edges changes degree parity only at $u$. Since the finite graph $K_u$ has an even number of odd-degree vertices, an odd number of its vertices must have odd degree in $G$. In particular, $K_u$ contains an odd-degree vertex of $G$. The same argument applies to $K_w$, giving two distinct odd-degree vertices of $G$, a contradiction. This shows that $K_u=K_w$ and therefore (3) can be satisfied, which is a contradiction.

%Observe that, since deleting the already oriented edges changes degree parity only at $u$ or $w$, if $K_u\neq K_w$ then each must contain an odd-degree vertex of $G$, against \thref{thm:EGW}. Hence $K_u=K_w$ and therefore (3) can be satisfied, which is a contradiction.

Since each stage $s$ of the construction is well-defined, every edge is eventually oriented. The construction yields an edge-decomposition of the resulting graph $D$ into a family $\{C_i\}_{i}$ of finite directed closed trails and at most one remaining directed trail $R$, which may be finite, one-way infinite, or two-way infinite. 

We now show that either $D$ or $D^\mathrm{op}$ admits a one-way directed Eulerian path. Observe that each closed directed trail $C_i$ contributes equally to the in-degree and out-degree of each vertex. The same holds for every vertex in $R$ different from the endpoint(s). At each endpoint, the contributions of $R$ differ by one, so any finite-degree endpoint has odd degree in $G$, while every infinite-degree vertex has infinite in-degree and out-degree.

To show that $D$ (and hence $D^\mathrm{op}$) is weakly unbiased, fix $X\subseteq V$. Each closed directed trail $C_i$ has equally many edges entering and leaving $X$, since either $C_i$ is contained in $X$ or $C_i$ is contained in $V\setminus X$, or every entrance in $X$ is followed by an exit. Likewise, along $R$, the edges entering and leaving $X$ alternate: indeed, assume $R$ is infinite (otherwise it does not affect weak unbiasedness) and one-way. If the tail of $R$ is contained in $X$ or in $V\setminus X$, then the set of edges of $R$ that enter and leave $X$ are finite. Otherwise, $X$ divides $R$ in infinitely many finite directed trails $R_i$, each one being a maximal consecutive subtrail of $R$ entirely contained in $X$ or in $V\setminus X$. Each $R_i\neq R_0$ has an incident edge entering $X$ and an incident edge leaving $X$ ($R_0$ is an exception as it only has an incident edge leaving $X$ or has an incident edge entering $X$), hence they do not affect the weak unbiasedness of $D$. The case when $R$ is a two-way directed trail is analogous. This shows that both $D$ and $D^\mathrm{op}$ are weakly unbiased.

%Thus infinitely many edges enter $X$ iff infinitely many leave it, which proves the unbiasedness. 

If $G$ has no odd-degree vertex, both orientations satisfy \thref{thm:nash-williams-one-way}, starting at any infinite-degree vertex. Otherwise, there is $v_0\in V$ such that $\deg_G(v_0)$ is odd. Then $\indeg_D(v_0)=\outdeg_{D^\mathrm{op}}(v_0)$ and $\outdeg_D(v_0)=\indeg_{D^\mathrm{op}}(v_0)$. In particular, exactly one of $D$ and $D^\mathrm{op}$ satisfies \thref{thm:nash-williams-one-way} starting at $v_0$. 
\end{proof}

\begin{theorem}
	\thlabel{thm:c2'}
	$\OrientOW\weireducible\Choice{2}'$.
\end{theorem}
\begin{proof}
	Given a graph $G$, apply \thref{thm:D_Dop} to compute two orientations $D_0$ and $D_1:=D_0^{\mathrm{op}}$ of $G$ such that at least one of them admits a one-way directed Eulerian path. By the proof of \thref{thm:D_Dop}, $D_i\notin\dom(\DEP)$ iff there is a vertex $v$ of finite degree such that $\indeg_{D_i}(v)=\outdeg_{D_i}(v)+1$. This is a $\Sigma^0_2$ condition relative to the input graph (see also the proof of \thref{thm:EP}). Hence the non-empty set
	\[ A:=\set{i\in\{0,1\}}{D_i \in \dom(\DEP)}\]
	is uniformly $\Pi^0_2$ relative to $G$, and $\Choice{2}'$ suffices to choose a valid orientation.
\end{proof}

As a corollary, we obtain e.g.\ that $\Choice{3}\not\weireducible\OrientOW$, since $\Choice{3}\not\weireducible\Choice{2}'$ \cite[Prop.\ 7.4]{BRramsey17}. We now show that $\OrientOW \weireducible \Sort_2$. The following observation is a simple consequence of \cite[Lem.\ 7]{NashWilliams60} (see also \cite[Thm.\ 8.4.2]{Diestel25}).

% We now show that $\OrientOW \weireducible \Sort_2$. Recall that $\Sort_2\function{\Cantor}{\Cantor}$ is defined as
% \[
% \Sort_2(p):=\begin{cases}0^n1^\omega&\text{if $p$ contains exactly $n$ zeroes},\\0^\omega&\text{if $p$ contains infinitely many zeroes}.\end{cases}
% \]
% It is known that  $\CNatural\strictlyweireducible\Sort_2\strictlyweireducible\mflim$ \cite[Prop.\ 24]{NP18}.

%The following observation is a simple consequence of \cite[Lem.\ 7]{NashWilliams60} (see also \cite[Thm.\ 8.4.2]{Diestel25}).

\begin{lemma}
	\thlabel{thm:edge-decomposition}
	Let $G\in \dom(\EP)$. If $G$ has no odd-degree vertex then it admits an edge-decomposition into finite cycles.
\end{lemma}
\begin{proof}
	By \cite[Thm.\ 8.4.2]{Diestel25}, it is enough to show that $G$ has no finite odd cut, i.e.\ there is no set $X\subseteq V$ such that the set $F$ of edges with exactly one endpoint in $X$ is finite and $|F|$ is odd. Fix $X$ and $F$ as mentioned, and suppose that $F$ is finite. Since $G$ is connected, removing $F$ leaves only finitely many connected components, and exactly one of them is infinite (\thref{thm:EGW}). Hence one side of the cut, say $X$, is finite. Every vertex in $X$ has finite even degree, hence, letting $E_{in}$ be the set of edges between two vertices in $X$, we have $\sum_{v\in X} \deg_G(v) = 2|E_{in}| + |F|$.
	Therefore $|F|$ is even, as required.
\end{proof}

The proof of \thref{thm:D_Dop} can be adjusted to show that, if we are given both a graph $G\in \dom(\EP)$ and a vertex $v_0$ such that $G$ admits a one-way Eulerian path starting from $v_0$, then we can computably choose which one between $D$ and $D^\mathrm{op}$ admits a directed Eulerian path. However, we will prove a more general lemma that will simplify the proof of $\OrientOW \weireducible \Sort_2$. 

We say that a function $\varphi\function{\mathbb{N}}{\mathbb{N}\cup \{-1\}}$ is an \textdef{enumerating function} for a set $A\subseteq \mathbb{N}$ if $\ran(\varphi) \setminus \{-1\} = A$. In other words, checking if $n\in A$ is a c.e.\ condition relative to $\varphi$. 

\begin{lemma}
	\thlabel{thm:orient+enum}
	There is a computable procedure that, given a graph $G=(V,E)\in \dom(\EP)$ and an enumerating function $\varphi$ for $\set{ v \in V}{\deg_G(v) \text{ is odd}}$, computes $D\in\OrientOW(G)$.
\end{lemma}

We first explain the strategy behind the construction. By \thref{thm:edge-decomposition}, every graph $G\in \dom(\EP)$ with no odd-degree vertex admits an edge-decomposition into finite cycles. A closer inspection of the proof readily reveals that the decomposition is effective. However, if there is an odd-degree vertex, a graph $G\in\dom(\EP)$ with an odd-degree vertex can only be edge-decomposed into $T \cup \bigcup_n C_n$, where each $C_n$ is a finite cycle and $T$ is a finite or one-way infinite trail. In this case, we use $\varphi$ to decide the orientation of $T$ so that the resulting directed graph is in $\dom(\DEP)$.

\begin{proof}
	Let $\sequence{e_s}{s\in\mathbb{N}}$ be an enumeration of $E$. At each stage $s$, we maintain a finite directed trail $\tau_s$, starting with $\tau_0:=\str{}$. If $e_s$ is already oriented by stage $s$, we just define $\tau_{s+1}:=\tau_s$ and go to the next stage. Let $H_s$ be the subgraph obtained from $G$ by removing the oriented edges and isolated vertices. We search in parallel for one of the following: 
	\begin{enumerate}
		\item a finite cycle $C$ in $H_s$ containing $e_s$;
		\item a stage $t\ge s$ such that $\varphi$ has enumerated a vertex $v_0$ by stage $t$. If such $t$ is found, we search for a finite trail $\rho$ starting at $v_0$ and containing $e_s$ such that $\rho$ extends $\tau_s$ and $E(\rho)\setminus E(\tau_s)\subseteq E(H_s)$.
	\end{enumerate}
	In case (1), we orient $C$ cyclically and go to the next stage with $\tau_{s+1}:=\tau_s$. In case (2), we orient $\rho$ pointing away from $v_0$ and go to the next stage with $\tau_{s+1}:=\rho$. 	
	
	We now show that the search always succeeds. If $G$ has no odd-degree vertex, then case (2) never occurs, and the claim follows by \thref{thm:edge-decomposition}. Assume instead that, at stage $s$, we find $t\ge s$ such that $v_0$ is enumerated by $\varphi$ by stage $t$. Let $w_s$ be the final vertex of $\tau_s$, taking $w_s:=v_0$ if $\tau_s$ is empty. Observe that every finite-degree vertex of $H_s$ different from $w_s$ has even degree, while $\deg_{H_s}(w_s)$ is odd if $\deg_G(w_s)<\infty$. Indeed, removing cycles preserves degree parity, and removing a trail from $v_0$ to $w_s$ changes parity only at its endpoints when these are distinct. If the endpoints coincide, or the trail is empty, the degree parity at $v_0$ is unchanged.

	Hence the connected component $K$ of $H_s$ containing $w_s$ is infinite: if $\deg_G(w_s)=\infty$, this is immediate. Otherwise, a finite $K$ would have exactly one odd-degree vertex, which is impossible. By \thref{thm:EGW}, $K$ is the unique infinite component of $H_s$. If $e_s$ belongs to $K$, there is a finite trail in $H_s$ starting at $w_s$ and containing $e_s$. Appending this trail to $\tau_s$ gives a $\rho$ satisfying (2). Otherwise, the component containing $e_s$ is finite with all degrees even, so $e_s$ belongs to a finite cycle and (1) applies. Thus the search succeeds in both cases.

	This shows that every edge is eventually oriented. We now prove that the resulting directed graph $D$ admits a one-way directed Eulerian path by showing that it satisfies the hypotheses of \thref{thm:nash-williams-one-way}. Clearly the underlying undirected graph of $D$ has an Eulerian path, and the same argument used in the proof of \thref{thm:D_Dop} shows that $D$ is weakly unbiased. If $G$ has no odd-degree vertex, the construction gives an edge-decomposition into finite directed cycles. In this case, every vertex $v$ satisfies $\indeg_D(v)=\outdeg_D(v)$. Otherwise, let $v_0$ be the unique odd-degree vertex of $G$. The construction gives an edge-decomposition into finite directed cycles and a finite or right-infinite directed trail $T$ starting at $v_0$. Observe that every directed trail contributes equally to the in-degree and out-degree of each vertex, except for its endpoint(s). If $T$ is right-infinite, this implies that the only unbalanced vertex is $v_0$ (for which $\indeg_D(v_0)+1=\outdeg_D(v_0)$). If $T$ is finite, we notice that its final vertex $w$ must have infinite degree. Indeed, suppose that $\deg_G(w)<\infty$. If $w\neq v_0$, then $T$ contributes an odd number of edges incident to $w$, while every closed trail contributes an even number. Hence $\deg_G(w)$ is odd, against the fact that $v_0$ is the only possible odd-degree vertex of $G$. If $w=v_0$ then $T$ is closed (possibly trivial), hence $\deg_G(v_0)$ is even, against the fact that $\deg_G(v_0)$ is odd. Thus the degree condition holds, completing the proof.
\end{proof}

We can notice that the construction in the previous Lemma can be adapted to show that, given $G\in\dom(\EP)$ and a vertex $v_0\in V(G)$ such that $G$ has an Eulerian path starting from $v_0$, we can uniformly compute a solution for $\OrientOW(G)$. Indeed, the same proof applies whether $v_0$ has odd or infinite degree. 

\begin{theorem}
	\thlabel{thm:sort}
	$\OrientOW\weireducible\Sort_2$.
\end{theorem}
\begin{proof}
	By \thref{thm:orient+enum}, it is enough to show that, given $G\in \dom(\EP)$, we can use $\Sort_2$ to compute an enumerating function for $\set{v \in V}{\deg_G(v) \text{ is odd}}$. Let $X$ be the set of pairs $(v,\sigma)\in V\times\baire$ such that $\length{\sigma}$ is odd and $\sigma$ lists exactly the neighbors of $v$, without repetitions. This set is uniformly $\Pi^0_1$ relative to $G$. Fix an enumeration of $V\times \baire$ and define $g(s)$ as the least pair that is not enumerated outside of $X$ by stage $s$. We define the input $p$ for $\Sort_2$ as 
	\[	p(s):=\begin{cases}
			0 & \text{if }g(s+1)\neq g(s),\\
			1 & \text{otherwise.}
		\end{cases}	\]
	% If $G$ has an odd-degree vertex, $g$ stabilizes at a pair $(v,\sigma)$ where $v$ has odd degree. Otherwise, $g$ changes infinitely often. Thus $q:=\Sort_2(p)$ is respectively $0^k1^\omega$, where $k$ counts the changes of $g$, or $0^\omega$. To compute $\varphi$, output $-1$ until the first $1$ in $q$ appears, then simulate $g$ until its $k$-th change (using $g(0)$ if $k=0$) and enumerate the vertex identified by its final value.
	%
	Observe that if $G$ has an odd-degree vertex, then $g$ eventually stabilizes at a pair $(v,\sigma)$ such that $v$ has odd degree and $\sigma$ is a list of the neighbors of $v$. Otherwise, every guess is eventually refuted, and $g(s)\to\infty$. In particular, $g$ changes its value finitely many times iff $G$ has an odd-degree vertex.	
	
	Let $q:=\Sort_2(p)$. Observe that if $G$ has no odd-degree vertex, then $p$ contains infinitely many zeroes and $q=0^\omega$. Otherwise, $q=0^k1^\omega$, where $k$ is the total number of changes of $g$.
	
	In particular, given $q$, we can compute an enumerating function $\varphi$ as follows: we output $-1$ as long as $q$ contains zeroes. Once the first $1$ in $q$ is found, using $k$ we can compute the final value of $g$, and hence the unique odd-degree vertex of $G$.
\end{proof}

\thref{thm:c2'} and \thref{thm:sort} immediately yield $\OrientOW \weireducible \Choice{2}' \sqcap \Sort_2$. We will now show that the reduction is strict. For the sake of readability, we first isolate the following lemma. 

\begin{definition}
	\thlabel{def:rotations}
Given a directed graph $F$, we can uniformly compute an enumeration $\sequence{C^F_i}{i}$ of all finite directed cycles in $F$. For each finite $\sigma\in\baire$, we define $\rho(\sigma, F)$ recursively as follows: $\rho(\str{},F):=F$ and $\rho(\str{i}, F)$ is the directed graph obtained by reversing the orientation of the $i$-th cycle $C^F_i$ in $F$, and $\rho(\sigma\concat \str{i}, F):= \rho(\str{i}, \rho(\sigma, F))$. 
\end{definition}

\begin{lemma}
	\thlabel{thm:finite_rotations}
	For any $G\in \dom(\EP)$, any $n\in \mathbb{N}$ and any two $D,H\in \OrientOW(G)$, there is $\sigma\in\baire$ such that $D_\sigma:=\rho(\sigma, D)$ is defined and agrees with $H$ on the orientation of the first $n$ edges of $G$. 
\end{lemma}
\begin{proof}
	Fix $G,D,H$ as in the hypotheses. Observe that, for every $\sigma\in\baire$ for which $D_\sigma$ is defined, $D_\sigma$ is a valid solution for $\OrientOW(G)$, as reversing the orientation of a cycle does not affect the existence of a one-way directed Eulerian path. We recursively define a string $\sigma$ as follows: we start with $\sigma_0:=\str{}$. At each stage, we have a set $A_s \subseteq E(D_s)$ of all the edges in $D_s$ with index $\le n$ on which $H$ and $D_s:=D_{\sigma_s}$ agree. Let $e=u\to v \in E(D_s)$ be the first edge with index $\le n$ on which $H$ and $D_{s}$ disagree (if no such $e$ exists we are done).

	Consider the set 
	\[  X:=\set{w\in V(G)}{\text{there is a directed path from }v \text{ to }w\text{ in }D_s-A_s}. \]
	If $u\in X$ then there is a cycle $C_i$ of $D_s$ proving the claim. We define $\sigma_{s+1}:=\sigma_s\concat \str{i}$ and go to the next stage. We show that $u\notin X$ leads to a contradiction. Indeed, for a directed graph $F$ and $X\subseteq V(F)$, let
	\begin{gather*}
		\delta_{\mathrm{out}}(F,X):=\set{(x,y)\in E(F)}{x \in X \text{ and }y\notin X} \\
		\delta_{\mathrm{in}}(F,X):=\set{(y,x)\in E(F)}{x \in X \text{ and }y\notin X}.
	\end{gather*}
	Clearly, $\delta_{\mathrm{out}}(D_s,X)\subseteq A_s$ hence, by definition of $A_s$, it is finite. Since $D_s$ is weakly unbiased, this implies also that $\delta_{\mathrm{in}}(D_s, X)$ is finite. Since $\delta(D_s, X):= \delta_{\mathrm{in}}(D_s, X) \cup \delta_{\mathrm{out}}(D_s, X)$ is finite, \thref{thm:EGW} implies that $G-\delta(D_s, X)$ (i.e.\ deleting the undirected edges in $\delta(D_s, X)$ from $G$) has exactly one infinite component.
	
	Let $S:= X$ if $X$ does not contain the infinite connected component of $G-\delta(D_s,X)$, and $S:=V(G)\setminus X$ otherwise. Since $\delta(D_s, X)=\delta(D_s, S)$, it follows that there are finitely many edges with an endpoint in $S$. In particular, every one-way directed Eulerian path eventually leaves $S$ permanently. Since both $D_s$ and $H$ are in $\dom(\DEP)$, we have that for $F\in \{D_s,H\}$,
	\[ |\delta_{\mathrm{out}}(F, S)| - |\delta_{\mathrm{in}}(F, S)| = \begin{cases}
	1 	& \text{if $S$ contains an odd-degree vertex of $G$}\\
	0 	& \text{otherwise.}
	\end{cases}
	\]
	This can be proved using e.g.\ \cite[Lem.\ 4]{NashWilliams66}. This, together with $|\delta(D_s, S)|=|\delta(H,S)|$ (which follows from the fact that $S$ has the same boundary in $G$), implies that $|\delta_{\mathrm{out}}(H,S)| = |\delta_{\mathrm{out}}(D_s,S)|$ and $|\delta_{\mathrm{in}}(H,S)| = |\delta_{\mathrm{in}}(D_s,S)|$. This, in turn, implies that  $|\delta_{\mathrm{out}}(H,X)| = |\delta_{\mathrm{out}}(D_s,X)|$. 
	
	Since $\delta_{\mathrm{out}}(D_s,X)\subseteq A_s$, it follows by definition of $A_s$ that $\delta_{\mathrm{out}}(D_s,X)\subseteq \delta_{\mathrm{out}}(H,X)$. However, by construction, $H$ and $D_s$ disagree on the direction of the edge $e$, hence $|\delta_{\mathrm{out}}(H,X)| \ge |\delta_{\mathrm{out}}(D_s,X)|+1$, which is a contradiction. 
\end{proof}

\begin{theorem}
	\thlabel{thm:orient<meet}
	$\OrientOW \strictlyweireducible \Choice{2}' \sqcap \Sort_2$.
\end{theorem}
\begin{proof}
	The reduction follows from $\OrientOW \weireducible \Choice{2}'$ (\thref{thm:c2'}), $\OrientOW\weireducible\Sort_2$ (\thref{thm:sort}), as $\sqcap$ is the meet in the Weihrauch lattice. 

	For the separation, suppose towards a contradiction that $\Choice{2}'\sqcap\Sort_2\weireducible\OrientOW$ via $\Phi,\Psi$. Recall that a solution for $\Choice{2}'\sqcap\Sort_2(p,q)$ is of the form $b\concat y$ where $b\in \{0,1\}$ determines whether $y\in \Choice{2}'(p)$ or $y=\Sort_2(q)$. It is convenient to use $\Choice{2}'\weiequiv \mathsf{RT}^1_2$ and therefore think of an instance of $\Choice{2}'$ as a binary sequence $p\in\Cantor$, and a $\Choice{2}'$ solution for $p$ is a number $<2$ that appears infinitely often in $p$. 
	
	For every instance $(p,q)$ of $\Choice{2}'\sqcap\Sort_2$, we write $G_{p,q}$ for the graph defined by $\Phi(p,q)$ and $D_0^{p,q}$, $D_1^{p,q}:=(D_0^{p,q})^\mathrm{op}$ for the two orientations of $G_{p,q}$ computably obtained by \thref{thm:D_Dop}. Consider the two computations $\Psi_i(p,q):=\Psi(p,q,D^{p,q}_i)$ for $i<2$. By \thref{thm:D_Dop}, at least one of the two computations produces a valid solution for $\Choice{2}'\sqcap\Sort_2(p,q)$. 
	
	We first claim that if, for every input $(p,q)$, every valid canonical orientation $D_i^{p,q}$ yields a $\Choice{2}'$-answer, then $\Choice{2}' \weireducible \OrientOW$. This would contradict \thref{thm:sort}, as $\Sort_2\weireducible \mflim$ \cite[Prop.\ 24]{NP18} but $\Choice{2}'\not\weireducible \mflim$ \cite[Thm.\ 11.7.54]{BGP17}.
	
	To prove the claim, assume that, for every input $(p,q)$ and $i<2$, if $D^{p,q}_i\in\OrientOW(\Phi(p,q))$ then $\Psi_i(p,q)(0)=0$. Given an instance $p$ of $\Choice{2}'$, produce the graph $G:=\Phi(p,1^\omega)$ and let $D\in\OrientOW(G)$. In parallel for every $\sigma\in\baire$ such that $D_\sigma:=\rho(\sigma,D)$ is defined (see \thref{def:rotations}), we run the computations $\Psi(p,1^\omega,D_\sigma)$ until we find some $\sigma$ such that $\Psi(p,1^\omega,D_\sigma)(0)\downarrow = 0$. Notice that the search is guaranteed to terminate by \thref{thm:finite_rotations}, because at least one canonical orientation $D_i^{p,1^\omega}$ is valid and, by assumption, its backward computation returns a $\Choice{2}'$-answer. Besides, every enumerated $D_\sigma$ is a valid solution for $\OrientOW(\Phi(p,1^\omega))$, hence this establishes $\Choice{2}'\weireducible\OrientOW$.

	This shows that we can choose an input $(p,q)$ and $i<2$ such that $\Psi_i(p,q)(0)\downarrow =1$. We now diagonalize against $\Phi, \Psi$ by showing that there is an input $(p_*, q_*)$ for $\Choice{2}'\sqcap \Sort_2$ on which the reduction fails to produce a valid answer. Let $n_0$ be sufficiently large so that $\Psi_i(p[n_0],q[n_0])(0)\downarrow=1$.
	
	We first show that there are two strings $\sigma\pprefix p$, $\tau\pprefix q[n_0]\concat 1^\omega$ such that, for some $k$, $10^k1 \prefix \Psi_i(\sigma,\tau)$. Indeed, 
	consider both computations $\Psi_j(p,q[n_0]\concat 1^\omega)$ for $j<2$. If there are no $\sigma,\tau$ as above, then $\Psi_i(p,q[n_0]\concat 1^\omega)\notin \Choice{2}'\sqcap \Sort_2(p,q[n_0]\concat 1^\omega)$. This implies that $\Psi_{1-i}(p,q[n_0]\concat 1^\omega)$ is the correct one instead. By continuity, there are sufficiently long prefixes $\sigma \pprefix p$ and $\tau \pprefix q[n_0]\concat 1^\omega$ on which $\Psi_{1-i}(\sigma,\tau)$ produces either $0b$ for some $b<2$ or $10^k1$ for some $k\in\mathbb{N}$. 
	In the first case, set $p_*:=\sigma\concat (1-b)^\omega$ and $\tau_*:=\tau$, otherwise set $p_*:=p$ and $\tau_*:=\tau\concat 0^{k+1}1$. In both cases, we are forcing $\Psi_{1-i}$ to produce the wrong $\Choice{2}'\sqcap \Sort_2$ answer on every input of the form $(p_*, q_*)$ where $q_*\extend \tau_*$. Hence, letting $h$ be the number of zeroes in $\tau_*$, there is $m\in\mathbb{N}$ such that  $10^h1 \prefix \Psi_i(p_*, \tau_*\concat 1^m)$. By continuity, this implies that the pair $\Phi,\Psi$ fails to produce the correct $\Choice{2}'\sqcap \Sort_2$ answer for $(p_*, q_*:=\tau_*\concat 1^m 0^\omega)$, against the definition of Weihrauch reducibility.
	
	To conclude the proof, fix $\sigma,\tau\in\cantor$ and $k\in\mathbb{N}$ as above. In this case, letting $q_1:=\tau \concat 0^\omega$, we obtain that for every $p\pextend \sigma$, $\Psi_i(p,q_1)$ does not produce a valid $\Choice{2}'\sqcap \Sort_2$ answer. 
	
	We distinguish two cases:
	\begin{itemize}
		\item if, for some $j$ and some sufficiently large $n_1$, $\str{0,j}\prefix \Psi_{1-i}(p[n_1],q_1[n_1])$, then, we let $p_*:=p[n_1]\concat (1-j)^\omega$ and $q_*:=q_1$.
		\item if, for some sufficiently large $n_1$, $\Psi_{1-i}(p[n_1],q_1[n_1])(0)=1$, then, by continuity, there are $t\in\mathbb{N}$ and $m\in\mathbb{N}$ such that $10^m1 \prefix \Psi_{1-i}(p,q_1[n_1]\concat 1^t)$ (such $m$ must exist because $q_1[n_1]\concat 1^\omega$ contains finitely many zeroes). We define $p_*:=p$ and $q_*:=q_1[n_1]\concat 1^t0^\omega$. 
	\end{itemize}
	In both cases, we obtain that neither $\Psi_0(p_*, q_*)$ nor $\Psi_1(p_*, q_*)$ are valid solutions, as required. 	
\end{proof}

\section{Two-way Eulerian paths}

We now turn our attention to the case of two-way Eulerian paths.

\begin{definition}
	Let $\TEP$ (resp.\ $\TDEP$) be the problem that takes as input an infinite countable undirected (resp.\ directed) graph that has a two-way Eulerian path, and outputs a sequence $\sequence{x_n}{n\in\mathbb{Z}}$ of vertices forming such a path. 
\end{definition}

Recall that, by \thref{thm:EP}, $\EP\weiequiv\DEP \weiequiv \CNatural'\times \mflim$. Intuitively, the first half ($\CNatural'$) comes from the difficulty of ``starting the path'', while the second half ($\mflim$) is the complexity of finding ``the tail of the path''. Obviously, there is no ``start'' in a two-way Eulerian path, and indeed we can show that $\EP_2\weiequiv \DEP_2 \weiequiv \mflim$. For the sake of readability, we split the proof in two results, one for the undirected and one for the directed case. 

\begin{theorem}
	\thlabel{thm:tep=lim}
	$\TEP\weiequiv \mflim$.
\end{theorem}
\begin{proof}
	The reduction $\mflim \weireducible \TEP$ can be proved with a simple modification of the graph used in the proof of $\mflim\weireducible \EPstart$ (\thref{thm:fixed-start}), where the only difference is that the vertices $\set{b_n}{n\in\mathbb{N}}$ are now arranged so as to form a two-way line $\hdots - b_{3} - b_{1} - b_{0} - b_{2} - b_{4} - \hdots$. 

	Let us now show that $\TEP\weireducible\mflim$. Let $G=(V,E)$ be an instance of $\TEP$. By \cite[Prop.\ 4.13]{CarrascoVargas25}, a finite trail $\tau$ with initial vertex $u$ and final vertex $v$ is extendible to a two-way Eulerian path iff
	\begin{enumerate}
		\item there are distinct edges $e,f\in E\setminus E(\tau)$ incident to $u$ and $v$, respectively;
		\item every connected component of $G-E(\tau)$ contains $u$ or $v$.
	\end{enumerate}
	We show that these conditions are uniformly $\Delta^0_2$ relative to $G$ and $\tau$. The first condition is $\Sigma^0_1$ relative to $G$ and $\tau$, since we can search for two distinct unused edges incident to the endpoints. The second condition is equivalent to requiring that, for every $w\in\ran(\tau)$, one of the following holds:
	\begin{enumerate}[label=(\Roman*)]
		\item $E(\tau)$ contains all the edges incident to $w$, or
		\item there is a finite trail from $w$ to $u$ or $v$ in $G-E(\tau)$.
	\end{enumerate}
	To see that, let $H:=G-E(\tau)$. Assume first that $(2)$ holds and let $w\in\ran(\tau)$. If $w$ has an incident edge that is not contained in $\tau$, then $w\in V(H)$, hence $(2)$ implies that there is a finite trail from $w$ to either $u$ or $v$. Conversely, fix a connected component $C$ of $H$. Since $G$ is connected, $C$ contains some vertex $w$ of $\ran(\tau)$: take any path from $C$ to $\ran(\tau)$. For this $w$, (I) fails (as otherwise $w$ would not be in $V(H)$), hence (II) holds, and therefore $C$ contains either $u$ or $v$. 

	Observe that the conditions (I) and (II) are, respectively, $\Pi^0_1$ and $\Sigma^0_1$ relative to $G$ and $\tau$. This implies that, using $\mflim \weiequiv\parallelization{ \LPO}$, we can check which finite trail $\tau$ is extendible to a two-way Eulerian path. 
	
	To conclude the proof, we observe that such a path can be computed (uniformly) by greedy search while ensuring that all the edges are visited: fix an enumeration $\sequence{e_n}{n\in\mathbb{N}}$ of $E$. Let $e_0:=\{u,v\}$. We start with $\tau_0(0):=u$ and $\tau_0(1):=v$. At stage $s+1$, we have defined a finite trail $\tau_s$ with endpoints $u_s$ and $v_s$. We search for a finite trail $\tau_{s+1}$ extending $\tau_s$ such that
	\begin{itemize}
		\item $\tau_{s+1}$ is extendible to an infinite two-way Eulerian path. 
		\item $e_{s+1}\in E(\tau_{s+1})$
		\item $\tau_{s+1}$ extends $\tau_s$ by adding at least one additional edge at each end.
	\end{itemize}
	By the inductive hypothesis, such a trail exists, and, in the limit, this construction yields a solution for $\TEP(G)$.
\end{proof}

To characterize the degree of $\TDEP$, we first prove an auxiliary result.

\begin{lemma}
	\thlabel{lem:directed-two-way-extension}
	Let $D=(V,E)$ be a directed graph admitting a two-way directed Eulerian path. Let $\tau$ be a finite directed trail and let  $u:=\tau(0)$ and $v:=\tau(\length{\tau}-1)$. Define the graph $H$ as the directed graph obtained from $D-E(\tau)$ by adding a fresh vertex $z$ and the edges $u\to z$ and $z \to v$.

	Then $\tau$ extends to a two-way Eulerian path in $D$ iff the underlying undirected graph of $H$ is connected.
\end{lemma}
\begin{proof}
	Suppose first that $\tau$ extends to a two-way directed Eulerian path in $D$. Replacing the occurrence of $\tau$ in this path by $u\to z\to v$ gives a two-way directed Eulerian path in $H$. In particular, the underlying undirected graph of $H$ is connected.

	Conversely, suppose that the underlying undirected graph of $H$ is connected. We show that $H$ admits a two-way directed Eulerian path. This implies that $\tau$ extends to a two-way Eulerian path in $D$: indeed, every two-way Eulerian path in $H$ must traverse the edges $u\to z \to v$. Replacing the finite trail $u\to z \to v$ with $\tau$ gives a two-way Eulerian path in $D$ extending $\tau$. 

	Let $X=\sequence{x_n}{n\in\mathbb{Z}}$ be a two-way Eulerian path in $D$ and let $\sigma$ be a finite subtrail of $X$ containing every edge of $\tau$. Define $a:=\sigma(0)$ and $b:=\sigma(\length{\sigma}-1)$. Without loss of generality, we can assume $a\neq b$ (extending $\sigma$ if necessary). Consider the finite directed subgraph $F$ of $H$ such that 
	\[ E(F)=(E(\sigma)\setminus E(\tau))\cup\{(u,z),(z,v)\}.\]
	Observe that, since the two trails $\tau$ and $u\to z \to v$ have the same endpoints, for every vertex $x\in V(F)\setminus\{a,b\}$
	\[ \outdeg_F(x)-\indeg_F(x) = 0. \]
	Besides, $\outdeg_F(a)=\indeg_F(a)+1$ and $\indeg_F(b)=\outdeg_F(b)+1$. This implies that $F$ can be decomposed into finitely many edge-disjoint cycles $C_0,\hdots, C_{k-1}$ and a finite open directed trail $\eta$ from $a$ to $b$. This can be proved e.g.\ applying \cite[Ex.\ 4.8]{B-JG09} to the graph obtained by adding a fresh auxiliary vertex $c$ and the edges $b\to c$ and $c\to a$ to $F$. 

	Let $Y_0=\sequence{y_n^0}{n\in\mathbb{Z}}$ be the two-way directed trail obtained by replacing $\sigma$ with $\eta$ in $X$. We extend $Y_0$ by attaching each cycle $C_i$ in stages as follows: at each stage $s<k$, we have a list of unused cycles. Since $H$ is connected there is $i<k$ such that, for some unused cycle $C_i=\sequence{w_n^i}{n< \length{C_i}}$ (with the assumption that $w_0^i=w_{\length{C_i}-1}^i$), $C_i$ has some vertex in common with $Y_s$. Choose $n$ such that $y^s_n$ belongs to $C_i$. We define $Y_{s+1}$ by visiting $C_i$ between $y^s_n$ and $y^s_{n+1}$. More formally, modulo reindexing $C_i$, we can assume $w_0^i=y^s_n$. We define
	\[ y^{s+1}_m:=\begin{cases}
		y^s_m  & \text{if }m\le n\\ 
		w_{m-n}^i  & \text{if }n < m < n+\length{C_i}\\ 
		y^s_{m-\length{C_i}+1} & \text{if } n+\length{C_i} \le m.
	\end{cases}\] 
	We then label $C_i$ as used and go to the next stage. Each stage is possible because $H$ is connected, and the procedure terminates after $k$ stages. The trail $Y_k$ is a two-way directed Eulerian path through $H$, which concludes the proof.
\end{proof}

\begin{theorem}
	\thlabel{thm:tdep=lim}
	$\TDEP\weiequiv \mflim$.
\end{theorem}
\begin{proof}
	The reduction $\mflim\weireducible \TDEP$ is obtained by a straightforward modification of the construction in \thref{thm:fixed-start} (see also the proof of \thref{thm:tep=lim}): arrange the vertices $\set{b_n}{n\in\mathbb{N}}$ in a two-way line $\hdots \to b_3 \to b_1 \to b_0 \to b_2 \to b_4 \to \hdots$. 

	For the converse reduction, let $D=(V,E)$ be an instance of $\TDEP$. For each nonempty finite directed trail $\tau$, let $H_\tau$ be the graph obtained by replacing $E(\tau)$ with $u\to z\to v$ as in \thref{lem:directed-two-way-extension}. We show that connectedness of the underlying undirected graph of $H_\tau$ is uniformly $\Delta^0_2$ relative to $D$ and $\tau$.

	Observe that $H_\tau$ is connected iff, for every $w\in\ran(\tau)$, at least one of the following holds:
	\begin{enumerate}
		\item every edge of $D$ incident to $w$ belongs to $E(\tau)$;
		\item there is a finite undirected path in $H_\tau$ from $w$ to $z$.
	\end{enumerate}
	For the forward implication, if (1) fails, then $w$ is a vertex of $H_\tau$, so connectedness gives (2). Conversely, let $G$ be the underlying undirected graph of $D-E(\tau)$. Since $D$ is connected, every connected component of $D-E(\tau)$ contains a vertex $w\in \ran(\tau)$. By construction, every such vertex $w$ has an edge that does not belong to $E(\tau)$. In particular, by $(2)$, $w$ is connected to $z$ in $H_\tau$. Since, clearly, $u$ and $v$ are connected to $z$, this shows that $H_\tau$ is connected.

	Conditions (1) and (2) are respectively $\Pi^0_1$ and $\Sigma^0_1$ relative to $D$ and $\tau$, hence, since $\ran(\tau)$ is finite, checking if $H_\tau$ is connected is uniformly $\Delta^0_2$ in $D$ and $\tau$. We can therefore use $\mflim\weiequiv\parallelization{\LPO}$ to check which $\tau$ extend to a two-way Eulerian path. 

	To complete the proof, fix an enumeration $\sequence{e_s}{s\in\mathbb{N}}$ of $E$. At each stage, we greedily search for an extension $\tau_{s+1}$ of the finite trail $\tau_s$ built so far such that 
	\begin{itemize}
		\item $\tau_{s+1}$ is extendible,
		\item $\tau_{s+1}$ contains $e_s$, 
		\item $\tau_{s+1}$ extends $\tau_{s}$ by adding at least one edge at each end.
	\end{itemize}
	We start with $\tau_0:=e_0$. Since $D\in \dom(\TDEP)$ the search is guaranteed to succeed at each stage. In the limit, we obtain a two-way trail passing through each edge, and therefore it is a valid solution for $\TDEP(D)$. 
\end{proof}

Similarly to the case of one-way Eulerian paths, we can consider the problem of orienting the edges of an undirected graph while preserving the existence of a two-way (directed or undirected) Eulerian path.

\begin{definition}
	Let $\OrientTW$ be the problem that, given $G=(V,E)\in\dom(\TEP)$, returns $D=(V,A)\in\dom(\TDEP)$ such that the map $\pi\function{A}{E}:=(v,w) \mapsto \{v, w\}$ is a (computable) bijection.
\end{definition}

Unlike what happens for the case of one-way Eulerian paths, finding an orientation of the edges of an undirected graph so that the resulting directed graph admits a two-way Eulerian path is computable.

\begin{theorem}
	\thlabel{thm:computable-orientation-two-way}
	The problem $\OrientTW$ is computable.
\end{theorem}
\begin{proof}
	Let $G\in\dom(\TEP)$. We claim that the algorithm used in \thref{thm:D_Dop} to build the directed graph $D$ gives a valid solution for $\OrientTW(G)$. Recall that the construction builds an edge-decomposition $D=R\cup\bigcup_{i\in I}C_i$, where each $C_i$ is a finite directed closed trail, $I\subseteq \mathbb{N}$ could be finite or infinite, and $R$ is a trail and could be empty, finite, one-way infinite, or two-way infinite. The construction maintains a finite directed open trail $\tau_s$ at each stage, and tries to progress the edge-decomposition by searching for a finite directed trail $\rho$ that only includes edges that have not been oriented so far such that $\rho$ satisfies one of the following ($e_n$ is the first edge that has not been oriented by stage $s$):
	\begin{enumerate}
		\item $\rho$ includes $e_n$ and is a closed trail edge-disjoint from $\tau_s$;
		\item $\rho$ includes $e_n$ and can be either prepended or appended to $\tau_s$ to form an open directed trail;
		\item $\rho$ can be appended to $\tau_s$ to form a closed trail.
	\end{enumerate}
	To show that the search always succeeds, let $H$ be the undirected graph of edges that have not been oriented by stage $s$, and let $K$ be the connected component of $H$ containing $e_n$. Let also $u$ and $v$ be the two endpoints of $\tau_s$. As in the proof of \thref{thm:D_Dop}, if neither (1), (2), nor (3) is possible, then $u\notin K$, $v\notin K$, and $K$ is infinite. Observe now that both $u$ and $v$ have unused incident edges, since the construction has removed an odd number of edges incident to each of them, while their degrees in $G$ are even or infinite. Assume that $u$ and $v$ belong to different connected components. The connected component $K_u$ of $u$ is infinite: if it were finite, $u$ would be the only odd-degree vertex in $K_u$, therefore the sum of the degrees of the vertices in $K_u$ would be odd, which is impossible for a finite graph. Likewise, the connected component $K_v$ of $v$ is infinite. Together with $K$, this would give three infinite components after deleting finitely many edges from $G$, contrary to \thref{thm:EGW2}. Thus $u,v$ lie in the same connected component of $H$, and (3) is possible.

	To conclude the proof, we need to show that $D$ is strongly unbiased. Fix $X\subseteq V(G)$. Let $\delta_{\mathrm{in}}(X)$ (resp.\ $\delta_{\mathrm{out}}(X)$) be the set of edges entering (resp.\ leaving $X$). Observe that, if we define, for a directed subgraph $F$,
	\begin{gather*}
		\delta_{F\to X}:= \set{(u,v)\in E(F)}{u\in V(F)\setminus X, v\in X }\\
		\delta_{X\to F}:= \set{(v,u)\in E(F)}{u\in V(F)\setminus X, v\in X },
	\end{gather*}
	then 
	\begin{gather*}
		|\delta_{\mathrm{in}}(X)| = |\delta_{R\to X}| + \sum_{i\in I}|\delta_{C_i \to X}|\\
		|\delta_{\mathrm{out}}(X)| = |\delta_{X\to R}| + \sum_{i\in I}|\delta_{X\to C_i}|.
	\end{gather*}
	Since each $C_i$ has equally many edges entering and leaving $X$, $\sum_{i\in I}|\delta_{C_i \to X}| = \sum_{i\in I}|\delta_{X\to C_i}|$. Consequently, if they are both infinite then both $|\delta_{\mathrm{in}}(X)|$ and $|\delta_{\mathrm{out}}(X)|$ are infinite. Assume now that they are both finite. Observe that, along a directed trail $R$, the edges entering and leaving $X$ alternate. Therefore, either $|\delta_{X\to R}|=|\delta_{R\to X}|=\infty$, and hence both $|\delta_{\mathrm{in}}(X)|$ and $|\delta_{\mathrm{out}}(X)|$ are infinite, or $||\delta_{\mathrm{out}}(X)|-|\delta_{\mathrm{in}}(X)||=||\delta_{X\to R}|-|\delta_{R\to X}||\le 1$. This shows that $D$ is strongly unbiased.

	Since $G$ (and hence the underlying undirected graph of $D$) admits a two-way Eulerian path, by \thref{thm:nash-williams-two-way} we conclude that $D$ admits a two-way directed Eulerian path, and hence $D\in\OrientTW(G)$.
\end{proof}

\bibliographystyle{mbibstyle}
\bibliography{bibliography}

\end{document}